\documentclass[pdflatex]{article}

\usepackage{amsmath}
\usepackage{amsfonts}
\usepackage{amssymb}
\usepackage{amsbsy}
\usepackage{amscd}
\usepackage{amsthm}
\usepackage{latexsym}
\usepackage{ulem}
\usepackage{algorithmic}
\usepackage[frenchb,english]{babel} %
\usepackage[latin1]{inputenc}
\usepackage[T1]{fontenc}
\usepackage[cyr]{aeguill} % pour les guillemets
\usepackage[pdftex]{graphicx}
\usepackage[all]{xy}
\usepackage{hhline}
\usepackage{exscale,cmmib57}
\usepackage{epic}%,eepic} % attention, dashline disparaît avec eepic
\usepackage{subfigure}

\graphicspath{
{/u/cermics/r/boyaval/Desktop/SPDE/CMAME/figures/},
{/u/cermics/r/boyaval/Desktop/homogenization/SIAM/figures/},
{/u/cermics/r/boyaval/Desktop/homogenization/homogenization_1D/figure/},
{/u/cermics/r/boyaval/Desktop/homogenization/homogenization_1D/old/},
{/u/cermics/r/boyaval/Desktop/stochRB/reducvar_controlvar/CMS/fig_pdf/}
}

\renewcommand{\emph}{\textit}
\renewcommand{\em}{\it}
\newcommand{\brk}[1]{\left(#1\right)}          % \brk{.}     => (.)
\renewcommand{\to}{\rightarrow}

\newcommand{\D}{\mathcal{D}}

\newcommand{\X}{{\boldsymbol{X}}}

\newcommand{\F}{{\boldsymbol{F}}}

\def\P{\mathbb{P}} %probabilit\'e, \'el\'ements finis
\def\F{\mathcal{F}} % tribu

\newcommand{\N}{{\mathcal N}}
\newcommand{\R}{\mathbb{R}}

\renewcommand{\div}{\operatorname{div}}

\newcommand{\grad}{\boldsymbol{\nabla}}

\DeclareMathOperator{\argmax}{argmax}

\newcommand{\Var}[1]{\mathbf{Var}\left(#1\right)}
\newcommand{\Span}[1]{\mathbf{Span}\left(#1\right)}
\usepackage{setspace}
\theoremstyle{plain}% default
\newtheorem{thm}{Theorem}[section]

\newtheorem{proposition}[thm]{Proposition}

\title{Reduced basis techniques for stochastic problems}

\author{
S\'{e}bastien Boyaval$^{\rm a,b}$,
Claude Le Bris$^{\rm c,b}$, T. Leli\`evre$^{\rm c,b}$, \\
Yvon Maday$^{\rm d,e}$, Ngoc Cuong Nguyen$^{\rm f}$ and Anthony T. Patera$^{\rm f}$ \\[2ex]
$^{\rm a}$
\normalsize{Universit\'e Paris Est, Laboratoire Saint-Venant}\\
\normalsize{(Ecole des Ponts ParisTech) 6 \& 8 avenue Blaise Pascal}\\
\normalsize{Cit\'e Descartes, 77455 Marne-la-Vall\'ee Cedex 2, France}\\[2ex]
$^{\rm b}$
\normalsize{INRIA, MICMAC team-project, Domaine de Voluceau, BP. 105 - Rocquencourt}\\
\normalsize{78153 Le Chesnay Cedex France}\\[2ex]
$^{\rm c}$
\normalsize{Universit\'e Paris Est, CERMICS}\\
\normalsize{(Ecole des Ponts ParisTech) 6 \& 8 avenue Blaise Pascal}\\
\normalsize{Cit\'e Descartes, 77455 Marne-la-Vall\'ee Cedex 2, France}\\[2ex]
$^{\rm d}$
\normalsize{UPMC Univ Paris 06, UMR 7598, Laboratoire Jacques-Louis Lions} \\ 
\normalsize{F-75005, Paris, France}\\[2ex]
$^{\rm e}$
\normalsize{Division of Applied Mathematics, Brown University}\\
\normalsize{Providence, RI 02912 USA} \\[2ex]
$^{\rm f}$
\normalsize{Massachusetts Institute of Technology, Dept.~of Mechanical Engineering,} \\
\normalsize{Cambridge, MA  02139 USA}\\[2ex]
}

\date{}

\begin{document}

\maketitle

\section*{Abstract}

We report  here on the recent application of a now classical general
reduction technique, the  \textit{Reduced-Basis} (RB) approach
initiated in~\cite{prud'homme02:_reliab_real_time_solut_param}, to the
specific context of  differential equations with random
coefficients.  After an elementary presentation of the approach, we
review two contributions of the authors: \cite{boyaval-lebris-maday-nguyen-patera-09}, which presents the
application of the RB approach for the discretization of a  simple second
order elliptic equation supplied with a random boundary condition, and \cite{boyaval-lelievre-09}, which
uses a RB type approach to reduce the variance in the Monte-Carlo
simulation of a stochastic differential equation.  We conclude the
review with some general comments and also discuss possible tracks for
further research in the direction.

\section{Introduction}

% {\bf [blabla stochastic etc, plusieurs appels, plusieurs resolutions d'edp.]}

% We report  here on the recent application of a now classical general
% reduction technique, the  \textit{Reduced-Basis} (RB) approach
% initiated in~\cite{prud'homme02:_reliab_real_time_solut_param}, to the
% specific context of  differential equations involving  random
% functions.

% {\bf The RB approach

% -many queries

% - successful in many contexts

% - the random context is typically many-query

% -thus natural extension
% }

In this work we describe reduced basis (RB) approximation and
{\em a posteriori\/} error estimation methods for rapid and
reliable evaluation of {\em input-output
relationships\/} in which the {\em output\/} is expressed as a
functional of a {\em field variable \/} that is the solution of
an {\em
input-parametrized\/} system. In this paper our emphasis is on stochastic phenomena: the parameter is random; the system is a partial differential equation with random coefficients, or a stochastic differential equation, namely a differential equation forced by a Brownian process.

The reduced basis approach is designed to serve two important, ubiquitous, and challenging engineering contexts: real-time, such as estimation and control; and many-query, such as design, multi-scale simulation, and --- our emphasis here --- statistical analysis.
The parametric real-time and many-query contexts represent not only
computational challenges, but also computational opportunities: we may restrict our attention to a manifold of solutions, which can be rather accurately represented by a {\em low-dimensional vector space}; we can accept greatly increased pre-processing or ``Offline'' cost in exchange for greatly decreased ``Online'' cost for each new input-output evaluation. (All of these terms, such as "Online," will be more precisely defined in the Section 2.1 which constitutes a pedagogical introduction to the reduced basis approach.) Most variants of the reduced basis approach exploit these opportunities in some important fashion.

Early work on the reduced basis method focused on deterministic
algebraic and differential systems arising in specific domains \cite{fox71:_approx_analy_techn_desig_calcu,almroth78:_autom,noor81:_recen,noor82,noor80:_reduc,noor83:_multip,noor83:_recen,nagy79:_modal}; the techniques were subsequently extended to more general finite-dimensional systems as well as certain classes of partial differential equations (and ordinary differential equations) \cite{barrett95:_reduc_basis_method,fink83,lee91:_estim,noor84:_reduc,noor84:_mixed,porsching87,rheinboldt81:_numer_analy_contin_method_nonlin_struc_probl,rheinboldt93:_theor_error_estim_reduc_basis,porsching85:_estim}; the next decades saw further expansion into different applications and classes of equations, such as fluid dynamics and the incompressible Navier-Stokes equations \cite{peterson89,gunzburger89:_finit}. There is ample evidence of potential and realized success.

Recent research in reduced basis methods for deterministic parametrized partial differential equations both borrows from earlier efforts and also emphasizes new components: sampling techniques for construction of optimal reduced basis approximation spaces in particular in higher dimensional parameter domains~\cite{sen-08,boyaval-lebris-maday-nguyen-patera-09,nguyen-07}; rigorous {\it a posteriori} error estimation in appropriate norms and for particular scalar outputs of interest~\cite{knezevic-patera-09,Haasdonk-08}; and fastidious separation between the offline stage and online stage of the computations to achieve very rapid response~\cite{nguyen-veroy-patera-05}. These reduced basis methods can now be applied to larger, more global parameter domains, with much greater certainty and error control.

In this paper we emphasize the application of certified reduced basis methods to stochastic problems. Two illustrative approaches are explored. In the first approach~\cite{boyaval-lebris-maday-nguyen-patera-09} we consider application of the reduced basis approach to partial differential equations with random coefficients: we associate realizations of the random solution field to deterministic solutions of a parametrized deterministic partial differential equation; we apply the classical reduced basis approach to the parametrized deterministic partial differential equation. Statistical information may finally be obtained, for example through Monte Carlo approximations. New issues arise related to the simultaneous approximation of both the input random field and the solution random field.

In the second approach~\cite{boyaval-lelievre-09} we directly consider a statistical embodiment of the reduced basis notions. Here reduced basis ideas originally conceived in the deterministic differential context are re-interpreted in the statistical context: the deterministic differential equation is replaced by a parametrized random process; snapshots on the parametric manifold are replaced by correlated ensembles on the parametric manifold; error minimization (in the Galerkin sense) is replaced by variance reduction; offline and online stages are effected through fine and coarse ensembles. This technique is here applied to parametrized stochastic differential equations.

We begin, in Section~\ref{sec:initiation-RB},  with an initiation to the RB approach, considering a
simple, prototypical  elliptic problem, with deterministic coefficients.
Section~\ref{sec:SPDE} then presents the application of the approach
to a boundary value problem supplied with a random boundary
condition. The section summarizes the results some of us obtained
in~\cite{boyaval-lebris-maday-nguyen-patera-09}. With
Section~\ref{sec:SDE}, we address a problem different in nature,
although also involving randomness. The issue considered is the variance
reduction of a Monte-Carlo method for solving a stochastic differential
equation. The RB approach has been successfully  employed in \cite{boyaval-lelievre-09} to  efficiently generate
companion variables that are used as control variate and eventually reduce the
variance of the original quantities. The section outlines the approach
and shows its success on representative results
obtained. We conclude the article presenting
in Section~\ref{sec:futur} some potential, alternate applications of the  approach in the random
context.

\section{An Initiation to Reduced-Basis Techniques}
\label{sec:initiation-RB}

We begin  with an overview of Reduced Basis
techniques. The level of our exposition is elementary. Our purpose here
is to introduce the main ideas underlying the approach, leaving aside
all unnecessary technicalities. The reader already familiar with this family of
approximation approaches may easily skip this section and directly
proceed to sections~\ref{sec:SPDE} and~\ref{sec:SDE} where the adaptation of the general technique
to the specific case of partial differential equations with random
coefficients and to variance reduction using the RB approach will be addressed. We also refer
to~\cite{Rozza08:arcme,quarteroni-09} for pedagogic introductions to the standard RB method,
though with different perspectives.

\medskip

\subsection{Outline of the Reduced Basis approach}
\label{ssec:outlineRB}

Assume  that we need to evaluate,  for many values
of the parameter $\mu$,  some {\it output} quantity
$s(\mu)=F(u(\mu))$   function of the solution $u(\mu)$ to a partial
differential equation parametrized by this parameter $\mu$. 
If the computation of $u(\mu)$ and $s(\mu)$
for each single value of the parameter $\mu$ already invokes elaborate
algorithms, and this is indeed the case in the context of partial
differential equations,
then the numerical simulation of $u(\mu)$ and $s(\mu)$ for many $\mu$
may become a computationally overwhelming task.
Reducing the cost of parametrized computations is thus a challenge
to the numerical simulation. This is the purpose of Reduced Basis
technique (abbreviated as RB throughout this article) to reduce this cost.

\medskip

Let us formalize our discussion in the simple case of a partial
differential equation which is an elliptic second order equation of the
form (see (\ref{eq:PB00}) below): 
$$
-\div\brk{ \uuline{A}(\mu) \grad u(\mu) } = f ,$$
on a domain $\D$ with homogeneous Dirichlet boundary conditions.
The mathematical setting is classical. We assume that the solution
of interest $u(\mu)\in X$ is an element of a Hilbert space $X$
with inner product $(\cdot,\cdot)_X$ and norm $\|\cdot\|_X$. The output
$s(\mu) = F(u(\mu))\in\R$ is a scalar quantity
where $F:X\to\R$ is a smooth  (typically linear) function
and $\mu$ is a $P$-dimensional parameter varying in a fixed
given range $\Lambda\subset\R^P$. An example of such output~$s$ is $\displaystyle s(\mu) = F(u(\mu)) := \int_{\D} f \: u(\mu) $ (see
(\ref{eq:outRBput00}) below). The 
function $u(\mu)$ is mathematically defined as the solution to the general variational formulation:
\begin{equation}
\label{eq:primal00}
\text{Find $u(\mu)\in X$ solution to }
a(u(\mu),v;\mu)=l(v)\,,\ \forall v\in X\,, 
\end{equation}
where $a(\cdot,\cdot;\mu)$ is a 
symmetric bilinear form, continuous and coercive on $X$
and where $l(\cdot)$ is a linear form, continuous on $X$.
For all $\mu\in\Lambda$, $a(\cdot,\cdot;\mu)$ thus defines an inner
product in $X$. The existence and uniqueness of $u(\mu)$, for each
$\mu$, is then obtained by standard arguments. 

We henceforth denote by $\|\cdot\|_\mu$ the  norm $\displaystyle
\|\cdot\|_\mu=\sqrt{a(\cdot,\cdot;\mu)}$ equivalent to $\|\cdot\|_X$ (under appropriate assumptions on $\uuline{A}$, see below),
which is usually termed the energy norm.
In the sequel, we denote by
$u_{\N}(\mu)\in X_\N$ an accurate Galerkin approximation for $u(\mu)$
in a linear subspace $X_\N \subset X$ of dimension $\N\gg1$
and by $s_{\N}(\mu)=F(u_{\N}(\mu))$ the corresponding approximation for the output $s(\mu)$.
For that particular choice of~$X_\N$,
we  assume that the approximation error $|s(\mu)-s_{\N}(\mu)|$
is uniformly sufficient small for all $\mu\in\Lambda$. That is,
 $s_{\N}(\mu)$ is considered as a good approximation of the output
 $s(\mu)$ in practical applications. The difficulty is, we put ourselves
 in the situation where computing $s_{\N}(\mu)$ for all the values $\mu$
 needed is too expensive, given the high dimensionality $\N$ of the
 space $X_\N$ and the number of parameters~$\mu$ for which equation~\eqref{eq:primal00} need to
 be solved.

The RB approach typically consists of two steps. The purpose of the
first step is to construct a 
linear subspace \begin{equation}\label{eq:defXN}
X_{\N,N} = \Span{u_\N(\mu^N_n),n=1,\ldots,N}\,,
\end{equation}
subset of $X_\N$, of dimension $N\ll\N$, using a few approximate
solutions to~\eqref{eq:primal00} for particular values of the parameter
$\mu$. The point is of course to carefully select these values
$(\mu^N_n)_{1\le n\le N}\in\Lambda^N$ of the parameter, and we will
discuss this below (see \eqref{eq:idealXN} and \eqref{eq:practicalXN}). For intuitively clear reasons, the particular
solutions $u_\N(\mu^N_n)$ are called \emph{snapshots}. This first step is called the \emph{offline}
step, and is typically an expensive computation, performed once for all. In a second step, called the \emph{online} step, an approximation
$u_{\N,N}(\mu) \in X_{\N,N}$ of the solution to~\eqref{eq:primal00}  
is computed as a linear combination of the  $u_\N(\mu^N_n)$. The
problem solved states:
\begin{equation}
\label{eq:primal00N}
\text{Find $u_{\N,N}(\mu)\in X_{\N,N}$ solution to }
a(u_{\N,N}(\mu),v;\mu)=l(v)\,,\ \forall v\in X_{\N,N}\,.
\end{equation}
This problem is much less computationally demanding than solving for the fine solution $u_{\N}(\mu)$, and will be performed for many values of the parameter $\mu$.
We denote by $s_{\N,N}(\mu)=F(u_{\N,N}(\mu))$ 
the corresponding approximation of the output $s(\mu)$. An 
 {\it a posteriori} estimator $\Delta_N^s(\mu)$ 
for the output approximation error $|s_{\N}(\mu)-s_{\N,N}(\mu)|$ is
needed in order to appropriately calibrate~$N$ and select the $(\mu^N_n)_{1 \le n \le N}$. This {\it a posteriori} estimator may also be used in the online step to check the accuracy of the output. We shall make this precise below. For the time being, we only emphasize that the {\em a posteriori} analysis we develop aims at assessing the quality of the approximation of the output $s$ (and not the overall quality of the approximation of the solution), see~\cite{ainsworth00:_poster_error_estim_finit_elemen_analy} and references therein. The method is typically called a {\em goal oriented} approximation method.

The formal argument that gives hope to construct an accurate
approximation of the solution $u(\mu)$ to \eqref{eq:primal00} using this
process is that the 
manifold $\mathcal{M}_{\N} = \{u_\N(\mu),\mu\in\Lambda\}$ is expected to
be well approximated by a linear space of dimension much smaller than $\N$, the dimension of the ambient space $X_\N$. An
expansion on \emph{a few} snapshots $N$ has therefore a chance to
succeed in accurately capturing the solution $u(\mu)$ for all parameter values~$\mu$. The reduced basis method is fundamentally a discretization method to appoximate the state space~$\mathcal{M}_{\N}$, with a view to computing an accurate approximation of the output. Of course, this requirement strongly depends on the choice of the parametrization which is a matter of modelling.

The RB method yields good approximations $s_{\N,N}(\mu)$ of $s_{\N}(\mu)$ under appropriate assumptions on the dependency of the solution $u(\mu)$ on the input parameter $\mu$.
As a consequence,
optimal choices for the approximation space $X_{\N,N}$ 
should  account for the dependency of the problem with respect to~$\mu$.
More precisely,
the method
should select parameter values $(\mu^N_n)_{1\le n\le N}\in\Lambda^N$
with a view to controlling the  norm of
the output approximation error $|s_{\N}(\mu)-s_{\N,N}(\mu)|$ as a function of $\mu$.
For most applications, the appropriate norm to consider for the error as a function of $\mu$ is the
$L^\infty$ norm  and this is the choice indeed made by the RB approach,
in contrast to many other, alternative approaches. The desirable  choice of
$(\mu^N_n)_{1\le n\le N}$ is thus defined by:
\begin{equation}
\label{eq:idealXN}
(\mu^N_n)_{1\le n\le N}
\in \underset{(\mu_n)_{1\le n\le N}\in\Lambda^N}{\mathop{\rm arginf}}
\left(\underset{\mu\in\Lambda}{\mathop{\rm sup}}\:|s_{\N}(\mu)-s_{\N,N}(\mu)|\right)
\end{equation}
Note that, although not explicitly stated,  the rightmost term
$s_{\N,N}(\mu)$ in~\eqref{eq:idealXN} parametrically depends on~$(\mu_n)_{1\le n\le N}$
because the solution to \eqref{eq:primal00} for $\mu$ is developped as a linear combination of the
corresponding snapshots~$u_{\N}(\mu_n)$.

It  is unfortunately very difficult to compute~\eqref{eq:idealXN} in practice.
With the publication~\cite{prud'homme02:_reliab_real_time_solut_param}, the RB approach
suggests an alternative, practically feasible procedure. Instead of
 the
 parameters $(\mu_n)_{1\le n\le N}$ defined by~\eqref{eq:idealXN}, the
 idea is to select   approximate minimizers of
\begin{equation}
\label{eq:practicalXN}
(\mu^N_n)_{1\le n\le N}
\in \underset{(\mu_n)_{1\le n\le N}\in\Lambda^N}{\mathop{\rm arginf}}
\left(\underset{\mu\in\Lambda_{\rm trial}}{\sup}\:\Delta_N^s(\mu)\right).
\end{equation}
Note that there are two differences between \eqref{eq:idealXN} and
\eqref{eq:practicalXN}. First, the set $\Lambda$ has been discretized into a very large trial sample of parameters 
$\Lambda_{\rm trial}\subset\Lambda$. Second, and more importantly, the
quantity $\Delta_N^s(\mu)$ minimized  in \eqref{eq:practicalXN} is an
\emph{estimator} of  $|s_{\N}(\mu)-s_{\N,N}(\mu)|$. A fundamental additional ingredient is that the approximate
minimizers of \eqref{eq:practicalXN} are selected using a specific
procedure, called \emph{greedy} because 
the parameter values $\mu^N_n$, $n=1,\ldots,N$, are selected
incrementally. Such an incremental procedure is in particular interesting
when $N$ is not known in advance,
since  the computation of approximate $\mu^N_n$ ($1\le n\le N$)  does not depend on $N$ 
and may be performed until the infimum in~\eqref{eq:practicalXN} is
judged sufficiently low.

Of course, the computation of approximations to~\eqref{eq:practicalXN} 
with such a greedy algorithm can still be expensive,
because a very large trial sample of parameters $\Lambda_{\rm trial}\subset\Lambda$
might have to  be explored. The RB method is thus only considered efficient when
 the original problem, problem~\eqref{eq:primal00} here,
has to be computed  for such a large number of input parameter values
$\mu$, that the overall procedure (computationally expensive offline step and then, efficient online step) is  practically more amenable
than following the
original, direct approach. One often speaks of a {\it many-query}
computational context when it is the case.  Notice that the RB method is
not to be seen as a competitor to the usual discretization methods;
it rather builds upon  already efficient  discretization methods 
using appropriate choices of $X_\N$ in order
to speed up   computations  that have to be performed repeatedly.

The approach can therefore be reformulated as the following two-step procedure 
\begin{itemize}
 \item in the {\it offline} stage (which, we recall, may possibly be computationally expensive),
one ``learns'' from a very large trial 
sample of parameters $\Lambda_{\rm trial}\subset\Lambda$
how to choose a small number $N$ of parameter values; this is performed
using a greedy algorithm that incrementally selects the $\mu_n$,
$n=1,\ldots,N$; the selection is based on the estimator~$\Delta_N^s(\mu)$;  accurate approximations $u_\N(\mu_n)$ for solutions $u(\mu_n)$ to~\eqref{eq:primal00} 
are correspondingly computed at those few parameter values;
 \item in the {\it online} stage,
computationally inexpensive approximations $u_{\N,N}(\mu)$ of solutions 
$u(\mu)$ to~\eqref{eq:primal00} are computed
for many values $\mu\in\Lambda$ of the parameter, using the Galerkin projection~\eqref{eq:primal00N}; the latter values need
not be in the  sample $\Lambda_{\rm trial}$,
and yield  approximations $s_{\N,N}(\mu)$ for the output $s(\mu)$; the
estimator  $\Delta_N^s(\mu)$, already useful in the offline step, is
again employed  to check the quality of the online approximation (this
check is called \emph{certification}).
\end{itemize}
Notice that the computation of the error estimator $\Delta_N^s(\mu)$ needs to be inexpensive, in order to be efficiently used on the very large trial sample in the offline stage, and for each new parameter values in the online stage.

One might ask why we proceed with the reduced basis approach rather than simply interpolate $s(\mu)$, given the few values $\{s_{\mathcal N}(\mu_1), \ldots, s_{\mathcal N}(\mu_N)\}$. There are several important reasons: first, we have rigorous error estimators based on the residual that are simply not possible based on direct interpolation; second,  these residuals and error bounds drive the greedy procedure; third, the state-space approach provides Galerkin projection as an "optimal" interpolant for the particular problem of interest; and fourth, in higher parameter dimensions (say of the order of 10 parameters), in fact the a priori construction of scattered-data interpolation points and procedures is very difficult, and the combination of the greedy and Galerkin is much more effective.

We are now in position to give some  more details on both the
offline and online steps of the RB approach in a very simple case: an elliptic problem, with an affine dependency on the parameter. Our next section will make specific what the
greedy algorithm, the estimator $\Delta_N^s(\mu)$, along with other objects
abstractly manipulated above, are.

\medskip

\subsection{Some more details on a simple case}
\label{ssec:outlineRB}

As mentioned above, we consider for simplicity the  Dirichlet
problem
\begin{equation}
 \label{eq:PB00}
\left\{ 
\begin{aligned} 
-\div\brk{ \uuline{A}(\mu) \grad u(\mu) } = f \text{ in } \D \,, %H^{-1}(\D)  \,,
\\
u(\mu) = 0 \text{ on } \partial\D \,, %H^{\frac12}(\partial\D) \,,
\end{aligned}
\right.
\end{equation}
where $\D$ is a two-, or three-dimensional domain
and the matrix $\uuline{A}(\mu)$ is parameterized by a single scalar
parameter $\mu\in\Lambda=[\mu_{min},\mu_{max}] \subset \R_+^*$. We assume that
the matrix $A$ is symmetric and depends on $\mu$ in an {\it affine} way:
\begin{equation}
\label{eq:affine00}
\uuline{A}(\mu) = \uuline{A_0} + \mu\:\uuline{A_1} \,, \ \forall \mu\in\Lambda \,.
\end{equation}
This assumption~\eqref{eq:affine00} is a crucial
ingredient, responsible, as we shall explain below, for a considerable
speed-up and thus  for the success of the RB approach here. More generally, either we must identify by inspection or construction an "affine" decomposition of the form~\eqref{eq:affine00}, or we must develop an appropriate affine approximation; both issues are discussed further below.

We assume we are interested in efficiently computing, for many values of $\mu\in\Lambda$, the output:
\begin{equation}
\label{eq:outRBput00}
s(\mu) = F(u(\mu)) := \int_{\D} f \: u(\mu)\,.
\end{equation}
This is of course only a specific situation. The output function can be
much more general, like a linear form $\displaystyle \int_{\D} g \,
u(\mu)$ with some $g\not=f$. Many other cases are possible, but they all
come at a cost, both in terms of analysis and in terms of workload. The
case \eqref{eq:outRBput00}, where the output coincides with the
linear form present in the right-hand side of the variational formulation of~\eqref{eq:PB00} (and where
the bilinear form $a$ involved in the variational formulation is symmetric), is called
\emph{compliant}. Having \eqref{eq:outRBput00} as an output function in
  particular simplifies the \emph{a posteriori} error analysis of the
  problem (namely the construction of $\Delta^s_N(\mu)$). 

We equip the problem, somewhat vaguely formulated in \eqref{eq:PB00}-\eqref{eq:affine00}-\eqref{eq:outRBput00} above, with the
appropriate mathematical setting that allow for all our necessary
manipulations  below to make sense. For consistency, we now briefly summarize this setting. The domain $\D$ is an open bounded connected domain with Lipschitz boundary $\partial\D$,
the right-hand side $f\in L^2(\D)$ belongs to the Lebesgue space of square
integrable functions, $\uuline{A}(\mu)$ is 
a  {\it symmetric} matrix, which is positive-definite almost everywhere
in $\D$. Each entry of $\uuline{A}(\mu)$ is assumed in $L^\infty(\D)$.
We  assume $\uuline{A_0}$ is symmetric positive-definite,
and $\uuline{A_1}$ is symmetric positive.
The ambient Hilbert space $X$ is chosen equal to  the usual Sobolev space $H^1_0(\D)$. The function  $u(\mu)$
is defined as the  solution to the variational formulation~\eqref{eq:primal00}
with $\displaystyle a(w,v;\mu)=\int_{\D}\uuline{A}(\mu)\nabla
w\cdot\nabla v$, $\displaystyle  l(v)=\int_{\D}f\,v$, for all $v$, $w$,
in $X$ and all $\mu\in\Lambda$. 

As for the direct discretization of the problem, we also put ourselves
in a classical situation. If $\D$ is polygonal for instance,
there exist many  discretization methods
that allow to compute Galerkin approximations $u_\N(\mu)$
of $u(\mu)$ in finite dimensional linear subspaces $X_\N$ of $X$
for any fixed parameter value $\mu\in\Lambda$. The Finite-Element
method~\cite{strang_fix,Ciarlet} is of course a good example.
Then, for each parameter value $\mu\in\Lambda$,
the numerical computation of $u_\N(\mu) = \sum_{n=1}^\N U_n(\mu) \phi_n
$ on the  Galerkin basis  $(\phi_n)_{1\le n\le\N}$ of $X_\N$ is achieved by solving a large  linear system
$$ \text{Find $U(\mu) \in\R^\N$ solution to } \uuline{B}(\mu) U(\mu) = b \,, $$
for the vector $U(\mu)=\left(U_n(\mu)\right)_{1\le n\le\N}\in\R^\N$,
where $b=\left(l(\phi_n)\right)_{1\le n\le\N}$ is a vector in $\R^\N$
and $\uuline{B}(\mu) = \uuline{B_0} + \mu\:\uuline{B_1}$ 
is a $\N\times\N$ real invertible matrix.
Note that the assumption of affine parametrization~\eqref{eq:affine00} makes possible, for
each parameter value $\mu$, the computation of 
the entries of the  matrix $\uuline{B}(\mu)$  in
$O(\N)$ operations (due to sparsity),  using the precomputed integrals
$({\uuline{B_q}})_{_{ij}}=\int_{\D}\uuline{A_q}\nabla \phi_i\cdot\nabla \phi_j $,
$i,j=1,\ldots,\N$ for $q=0,1$.
The evaluation of $U(\mu)$ for many $J\gg 1$ parameter values $\mu$
using iterative solvers
costs $J\times O(\N^k)$ operations 
with $k \le 3$~\cite{golub-vanloan-96}, where $k$ depends on the sparsity and the conditioning number of the involved matrices.

\medskip

As mentioned above in our general, formal presentation, the goal of the
RB approach is  to build a smaller 
finite dimensional approximation space $X_{\N,N}\subset X_\N$
sufficiently good for all $\mu\in\Lambda$, with $N\ll\N$,
so that the computational cost is approximately reduced to 
$N\times O(\N^k) + J\times O(N^3)$, 
where $N\times O(\N^k)$ is the cost of offline computations
and $J\times O(N^3)$, the cost of online computations,
is {\it independent of $\N$}, using the Galerkin approximation~\eqref{eq:primal00N} in $X_{\N,N}$.

\medskip

We now successively describe in the following three paragraphs the
construction of the \emph{a posteriori} estimator, that of the 
\emph{greedy} algorithm employed in the offline step, and the combination of
all ingredients in the online step.

\subsubsection{A posteriori estimator.}

For the coercive elliptic problem~\eqref{eq:PB00}, the {\it a
  posteriori} error estimator $\Delta_N^s(\mu)$
for the output RB approximation error $|s_\N(\mu)-s_{\N,N}(\mu)|$ is  simple to devise,
based on a global {\it a posteriori} error estimator $\Delta_N(\mu)$
for $\|u_\N(\mu)-u_{\N,N}(\mu)\|_\mu$ using a classical technique with residuals.

We refer
to~\cite{veroy_lions,veroy03:_poster_error_bound_reduc_basis,Patera_Huynh06,deparis07,Patera_Ronquist07,nguyen-09,Calcolo}
for the construction of similar {\it a posteriori} error estimators
in various applied settings of the RB method.

\medskip

We first define the residual bilinear form 
$$ g(w,v;\mu)=a(w,v;\mu)-l(v)\,,\ \forall w,v\in X\,,\ \forall \mu\in\Lambda\,, $$
and  the operator $G(\mu):X_\N\to X_\N$ such that
$$ g(w,v;\mu)=\left(G(\mu)\:w,v\right)_X\,,\ \forall w,v\in X_\N\,,\
\forall \mu\in\Lambda\,. $$

We next assume we are given, for all $\mu\in\Lambda$,
a lower bound $\alpha_{LB}(\mu)$ 
for the coercivity constant of $a(\cdot,\cdot;\mu)$ on $X_{\cal N}$, that is,
\begin{equation}
\label{eq:infinf}
0 < \alpha_{LB}(\mu) \le \alpha_c(\mu) = \underset{w\in X_{\cal N}\backslash\{0\}}{\inf}\frac{a(w,w;\mu)}{\|w\|_X^2} \,,\ 
\forall \mu\in\Lambda\,.
\end{equation}
The lower bound $\alpha_{LB}(\mu)$ can be given by an {\it a priori} analysis before discretization
($\alpha_{LB}(\mu)$ would then be the coercivity constant of $a(\cdot,\cdot;\mu)$ on $X$),
or numerically evaluated based on an approximation procedure, which might be difficult in some cases, see~\cite{huynh-rozza-sen-patera-07,Rozza08:arcme}.

Then the a posteriori estimator we use is defined in the following.
\begin{proposition}
For any linear subspace $X_{\N,N}$ of $X_\N$,
there exists a computable error bound $\Delta^s_N(\mu)$ such that:
\begin{equation}\label{eq:bound11}
|s_\N(\mu)-s_{\N,N}(\mu)|\le 
\Delta^s_N(\mu) := \frac{\|G(\mu)\:u_{\N,N}(\mu)\|_X^2}{\alpha_{LB}(\mu)}\,,\ 
\forall \mu\in\Lambda \,.
\end{equation}
\end{proposition}

For consistency, we now briefly outline the  proof of this
proposition. We simply observe the sequence of equalities
\begin{align}
\nonumber
|s_{\N}(\mu)-s_{\N,N}(\mu)| & = 
|F(u_\N(\mu))-F(u_{\N,N}(\mu))| 
\\
\nonumber
& = |l(u_\N(\mu)-u_{\N,N}(\mu))|
\\
\nonumber
& = |a(u_\N(\mu),u_\N(\mu)-u_{\N,N}(\mu);\mu)|
\\
\nonumber
& = |a(u_\N(\mu)-u_{\N,N}(\mu),u_\N(\mu)-u_{\N,N}(\mu);\mu)|
\\
&=\|u_\N(\mu)-u_{\N,N}(\mu)\|_{\mu}^2\,.
\label{eq:estimate01}
\end{align}
using the linearity of $F=l$, 
the  variational problem and its discretized approximation in $X_\N$,
the symmetry of $a(\cdot,\cdot;\mu)$ and the fact that
$a(u_\N(\mu)-u_{\N,N}(\mu),v)=0$, for all $v\in X_{\N,N}$. 
On the other hand, inserting $v=u_\N(\mu)-u_{\N,N}(\mu)$ in the general
equality 
$a(u_\N(\mu)-u_{\N,N}(\mu),v;\mu)
=-g(u_{\N,N}(\mu),v;\mu)$ (for all $v\in X_\N$), and using the bound
$\sqrt{\alpha_{LB}(\mu)}\|v\|_X\le\|v\|_\mu$ (for all $v\in X_\N$), we
note that
\begin{equation}\label{eq:bound01}
\|u_\N(\mu)-u_{\N,N}(\mu)\|_{\mu} \le \Delta_N(\mu) := \frac{\|G(\mu)\:u_{\N,N}(\mu)\|_X}{\sqrt{\alpha_{LB}(\mu)}}\,.
\end{equation}
We conclude the proof of \eqref{eq:bound11} combining~\eqref{eq:estimate01} with~\eqref{eq:bound01}.

\medskip

We may similarly prove (but we will omit the argument here for brevity)  the  inverse inequality:
\begin{equation}\label{eq:sharpness}
\Delta^s_N(\mu) \le 
\left(\frac{\gamma(\mu)}{\alpha_{LB}(\mu)}\right)^2
|s_{\N}(\mu)-s_{\N,N}(\mu)|\,,
\end{equation}
using the continuity constant 
\begin{equation}
\label{eq:supsup}
\gamma(\mu) = \underset{w\in X_\N\backslash\{0\}}{\sup}
\underset{v\in X_\N\backslash\{0\}}{\sup} \frac{a(w,v;\mu)}{\|w\|_X \|v\|_X} \,
\end{equation}
of the bilinear form $a(\cdot,\cdot;\mu)$ on $X_\N$ for all~$\mu\in\Lambda$,
which is bounded above by the continuity constant on $X$.
The inequality~\eqref{eq:sharpness}
ensures sharpness of the {\it a posteriori} estimator~\eqref{eq:bound11},
depending of course on the quality of the lower-bound $\alpha_{LB}(\mu)$.

\subsubsection{Offline stage and greedy algorithm.}\label{sec:offline}

The greedy algorithm employed to select the snapshots~$u_{\N}(\mu_n)$ typically reads:
\begin{algorithmic}[1]
\STATE choose $\mu_1\in\Lambda$ randomly
\STATE compute $u_\N(\mu_1)$ to define $X_{\N,1}=\Span{u_\N(\mu_1)}$
\FOR {$n=2$ to $N$}
	\item choose 
$\mu_{n} \in \mathop{\rm argmax}\{\Delta_{n-1}^s(\mu),\mu\in\Lambda_{\rm trial}\}$
	\item compute $u_\N(\mu_n)$ to define $X_{\N,n} = \Span{u_\N(\mu_m)\,,\ m=1,\ldots,n}$
\ENDFOR
\end{algorithmic}
In the initialization step, we may equally use 	$\mu_1\in\mathop{\rm
  argmax}\{|s(\mu)|\,,\ {\mu\in\Lambda_{\rm small trial}}\}$, 
where $\Lambda_{\rm small trial}\subset\Lambda$ is  a very small trial
sample in $\Lambda$, much smaller than $\Lambda$ itself. Likewise, 
the algorithm can  in practice be terminated when the output approximation error is
judged sufficiently small (say, $|\Delta_N^s(\mu)|\leq\varepsilon$ for
all $\mu\in\Lambda_{\rm trial}$), and not when the iteration number reaches a
maximum  $n=N$. 

The choice of the trial sample $\Lambda_{\rm trial}$
(and similarly, the smaller sample $\Lambda_{\rm small trial}$) is a delicate
 practical issue. It is often simply taken as a random sample in $\Lambda$. 
Of course, this first guess may be insufficient
to reach the required accuracy level $\varepsilon$ in $\Delta_N^s(\mu)$, for all $\mu \in \Lambda$, in the online stage.
But fortunately, if the computation of $\Delta_N^s(\mu)$ for any
$\mu\in\Lambda$ is sufficiently inexpensive,
one can check this accuracy online for each query in $\mu$.
Should $\Delta_N^s(\mu)>\varepsilon$ occur for some online value of the parameter~$\mu$,
one can still explicitly compute $u_\N(\mu)$ for that exact same $\mu$
and enrich the space~$X_{\N,N}$ correspondingly. This bootstrap approach of course allows to  reach the required accuracy level~$\varepsilon$ at that~$\mu$. It provides significant computational reductions in the online stage
\emph{provided that}
the RB approximation space $X_{\N,N}$ does not need to be enriched too often online.
We will explain the methodology for fast computations of $\Delta_N^s(\mu)$ below.

The offline selection procedure needs to be consistent with the online
procedure, and thus the above greedy algorithm uses the same estimator
$\Delta_N^s(\mu)$ for all $\mu\in\Lambda_{\rm trial}$ as the online procedure.
Since the computation of $\Delta_N^s(\mu)$ is, by construction and on
purpose,  fast for all 
$\mu\in\Lambda$,
the exploration of a very large training sample 
$\Lambda_{\rm trial}$ (which is a subset of $\Lambda$) is possible offline.

No systematic procedure seems to be available, which allows to build good
initial guesses $\Lambda_{\rm trial}$ \emph{ex nihilo}.
Even for a specific problem, we are not aware either of any {\it a priori} results
that quantify how good an initial guess $\Lambda_{\rm trial}$ is .
The only option is, as is indeed performed by the RB approach, to {\it a posteriori} check, 
and possibly improve, the quality of the initial guess $\Lambda_{\rm trial}$
(however, the quality of the initial guess $\Lambda_{\rm trial}$ can be slightly improved
offline by using adaptive training samples in the greedy algorithm~\cite{Haasdonk-Ohlberger-08-2}).

The estimators $\Delta_N^s(\mu)$ are
employed in the greedy algorithm to filter candidate values for~$\Lambda$. Numerous numerical evidences support the success of this pragmatic
approach~\cite{veroy_lions,veroy03:_poster_error_bound_reduc_basis,Patera_Huynh06,patera_rozza,Patera_Ronquist07,deparis07,Calcolo,nguyen-09}.

Last, notice that the cost of offline computations scales as 
$W_{\rm offline} = O(|\Lambda_{\rm trial}|)\times \left(\sum_{n=1}^{N-1} w_{\rm online}(n)\right) 
+ N\times O(\N^k)$
where $w_{\rm online}(n)$ is the marginal cost of one online-type computation for
$u_{\N,n}(\mu)$ and $\Delta^s_n(\mu)$ at a selected parameter value $\mu\in\Lambda_{\rm trial}$
(where $1\le n\le N-1$),
and $O(|\Lambda_{\rm trial}|)$ includes a $\max$-search in $\Lambda_{\rm trial}$.
(Recall that $k\le 3$ depends on the solver used for large sparse linear systems.)

\subsubsection{Online stage: fast computations including {\em a posteriori} estimators.}

We now explain how to efficiently compute $u_{\N,n}(\mu)$, $s_{\N,n}(\mu)$ and $\Delta^s_n(\mu)$ 
once the RB approximation space $X_{\N,n}$ has been constructed.
This task has to be completed twice in the RB approach. First, this is used in 
the many offline computations when $\mu\in\Lambda_{\rm trial}$
explores the trial sample
in order to find $\mu_{n} \in \mathop{\rm argmax}\{\Delta_{n-1}^s(\mu),\mu\in\Lambda_{\rm trial}\}$
at each iteration $n$ of the greedy algorithm. Second, this
is used for the many online computations (when $n=N$).
We present the procedure in the latter case.

By construction, the family $\left(u_{\N}(\mu_n)\right)_{1\le n\le N}$ generated by the greedy algorithm described in Section~\ref{sec:offline} is a basis of $$X_{\N,N} = \Span{u_{\N}(\mu_n)\,,\ n=1,\ldots,N}.$$
For any $\mu\in\Lambda$,
we would then like to compute the RB approximation
$ u_{\N,N}(\mu) = \sum_{n=1}^N U_{N,n}(\mu) u_{\N}(\mu_n) $,
which can be achieved by solving a small $N\times N$ (full) linear system
$$\uuline{C}(\mu) U_N(\mu) = c,$$ for the vector
$U_N(\mu)=\left(U_{N,n}(\mu)\right)_{1\le n\le N}\in\R^N$, with  $c=\left(l(u_{\N}(\mu_n))\right)_{1\le n\le N}$  a vector in $\R^N$
and $\uuline{C}(\mu) = \uuline{C_0} + \mu\:\uuline{C_1}$ 
is a $N\times N$ real invertible matrix.
In practice, the matrix $\uuline{C}(\mu)$ is close to a singular matrix,
and it is essential to  compute the RB approximation as
$ u_{\N,N}(\mu) = \sum_{n=1}^N \tilde U_{N,n}(\mu) \zeta_n $
using a basis  $(\zeta_n)_{1\le n\le N}$  of $X_{\N,N}$
that is orthonormal for the inner-product $(\cdot,\cdot)_X$. The
determination of appropriate $(\zeta_n)_{1\le n\le N}$  is
easily performed, since $N$ is small, using 
Simple or Modified Gram-Schmidt procedures. The problem to solve states:
\begin{equation}
\label{eq:linearsystem}
\text{Find $\tilde U_N(\mu) \in\R^N$ solution to } 
\uuline{\tilde C}(\mu) \tilde U_N(\mu) = \tilde c \,,
 \end{equation}
where $\tilde U_N(\mu)=\left(\tilde U_{N,n}(\mu)\right)_{1\le n\le
  N}\in\R^N$, 
$\tilde c=\left(l(\zeta_n)\right)_{1\le n\le N}$ is a vector in $\R^N$
and $\uuline{\tilde C}(\mu) = \uuline{\tilde C_0} + \mu\:\uuline{\tilde C_1}$
is a $N\times N$ real invertible matrix.
So, for each parameter value $\mu\in\Lambda$,
the entries of the latter matrix $\uuline{\tilde C}(\mu)$ can be computed in
$O(N^2)$ operations using the precomputed integrals
$(\uuline{\tilde C_q})_{_{ij}}=\int_{\D}\uuline{A_q}\nabla \zeta_i\cdot\nabla \zeta_j $,
$i,j=1,\ldots,N$ for $q=0,1$.
(Note that the assumption of affine parametrization is essential here).
And the evaluation of $\tilde U_N(\mu)$ for many $J\gg 1$ parameter values $\mu\in\Lambda$
finally costs $J\times O(N^3)$ operations %(multiplications and additions).
using exact solvers for symmetric problems like Cholesky~\cite{golub-vanloan-96}.

For each $\mu\in\Lambda$, the output $s_{\N,N}(\mu)=F(u_{\N,N}(\mu))$
can also be computed very fast in $O(N)$ operations
upon noting that $F$ is linear and all the values $F(u_{\N,N}(\mu_n))$, $n=1,\ldots,N$ can be precomputed offline.
The corresponding {\it a posteriori} estimator $\Delta^s_N(\mu)$ given
by~\eqref{eq:bound11} has now to be computed, hopefully equally fast.
Because of  the affine dependence of $\uuline{A}(\mu)$ on $\mu$,
a similar affine dependence $G(\mu)=G_0 +\mu\:G_1 $ holds for the operator $G$ (for all $\mu \in \Lambda$), 
where $\displaystyle (G_0\:w,v)_X=\int_{\D}\uuline{A_0}\nabla
w\cdot\nabla v-\int_{\D}f\,v$, and $\displaystyle
(G_1\:w,v)_X=\int_{\D}\uuline{A_1}\nabla w\cdot\nabla v$ for all $v$,
$w$, in $X_\N$. 
So one can evaluate very fast the  norm
\begin{align}
\|G(\mu)\:u_{\N,N}(\mu)\|_X^2 &=
\|G_0\:u_{\N,N}(\mu)\|_X^2
+ 2\mu (G_0\:u_{\N,N}(\mu),G_1\:u_{\N,N}(\mu))_X \nonumber \\
& \quad
+ \mu^2 \|G_1\:u_{\N,N}(\mu)\|_X^2\label{eq:quadraticmu}
\end{align}
for $\mu\in\Lambda$, 
once, with obvious notation, the scalar products $(G_i\:u_{\N,N}(\mu_p),G_j\:u_{\N,N}(\mu_q))_X$, have been precomputed offline and stored.
Assuming that the lower-bound $\alpha_{\rm LB}(\mu)$ used in~\eqref{eq:bound11}
is known,
the computation of the {\it a posteriori} estimator $\Delta^s_N(\mu)$ itself is thus also very fast.
Notice that the {\it affine} parametrization~\eqref{eq:affine00} plays a crucial role in the
above decomposition of the computation.

Finally, the marginal cost of one online-type computation on $X_{\N,n}$ for one parameter value $\mu$
is $w_{\rm online}(n)=O(n^3)$ (where $n=1,\ldots,N$).
So, assuming that no basis enrichment is necessary during the online stage using
the RB approximation space $X_{\N,N}$
(that is, $\Delta_N^s(\mu)<\varepsilon$ for all the parameter values $\mu$ queried online),
the total online cost for many $J\gg 1$ parameter values $\mu$
scales as $W_{\rm online}=J\times O(N^3)$.
And, the total cost of computations with the RB approach is then
$W_{\rm offline} + W_{\rm online} 
= N\times O(\N^k) + \left(J+O(|\Lambda_{\rm trial}|)\right)\times O(N^3)$,
which has to be compared to $J\times O(\N^k)$ operations for a direct approach
(with $k\le 3$ depending on the solver used for large sparse linear systems). In the limit of infinitely many online evaluations $J \gg 1$, the computational saving of the RB approach is tremendous.

\subsection{Some elements of analysis, and some extensions, of the RB method}
\label{ssec:theory}

\subsubsection{Some elements of theory.}

The RB approach has undoubtedly proved successful in a large variety of
applications~\cite{veroy_lions,veroy03:_poster_error_bound_reduc_basis,Patera_Huynh06,patera_rozza,Patera_Ronquist07,deparis07,Calcolo,nguyen-09}.
The theoretical understanding of the approach is however still limited,
and is far from covering all practical situations of interest. Of
course, little theory is to be expected in the usual {\it a priori}
way. As  already explained,
the RB approach is deliberately designed to  {\it a posteriori}
adapt to practical settings. The only available {\it a priori} analysis is related to two
issues: the expected ``theoretical'' quality of the RB approximation, and the efficiency
of the greedy algorithm. We now briefly summarize what is known to date
on both issues. 

The RB approach is in fact expected to perform \emph{ideally}, in the following
sense. In the context of our simple problem~\eqref{eq:PB00}, it is
possible, adapting the classical Lagrange interpolation theory 
to the context of parameterized boundary value problems and assuming
that the matrix $\uuline{A_1}$ is non-negative,
to obtain an upper bound of~\eqref{eq:idealXN}. The following
theoretical {\it a priori} analysis result
follows. It states the exponential accuracy of the RB
approximation in terms on the dimension $N$ of the reduced basis.
\begin{proposition}%[{\it A priori} existence]
\label{prop:apriori_rbRB}
For all parameter ranges $\Lambda:=[\mu_{\min},\mu_{\max}]\subset\R_+^*$,
there exists an integer $N_0 = O\left(\ln\left(\frac{\mu_{\max}}{\mu_{\min}}\right)\right)$
as $\frac{\mu_{\max}}{\mu_{\min}}\to+\infty$,
and a constant $c>0$ independent of $\Lambda$ such that,
for all $N\ge N_0\ge 2$, there exist % a sample $ \Lambda_N $ of 
$N$ parameter values 
$ \mu_{\min} =: \lambda_1^N < \ldots < \lambda_n^N < \lambda_{n+1}^N < \ldots < \lambda_N^N:=\mu_{\max} $,
$n=2,\ldots,N-2$, sastisfying 
(recall $\| \cdot \|_{0} = \| \cdot \|_{\mu}$ with $\mu=0$ is an Hilbertian norm on $X$):
\begin{align}
\underset{\mu\in\Lambda}{\sup} &\left(
\inf\left\{ 
\| u_\N(\mu)-w \|_{0}\,,\,w\in \Span{u_\N(\lambda_n^N),n=1,\ldots,N}
\right\}
\right) \nonumber \\
& \le e^{-\frac{c}{N_0-1}(N-1)}
\sup_{\mu\in\Lambda} \| u_\N(\mu) \|_{0} \,.\label{eq:decayboundRB}
\end{align}
\end{proposition}
We refer to~\cite[Chapter 4]{boyaval-these} and~\cite{{patera_rozza}} for the proof of
Proposition~\ref{prop:apriori_rbRB}.

\medskip

The approximation space $\Span{u_\N(\lambda_n^N),n=1,\ldots,N}$ 
used  for the statement and the proof of Proposition~\ref{prop:apriori_rbRB} 
is different from the RB approximation space $X_{\N,N}$ 
built in practice by the RB greedy algorithm.
Numerical experiments even suggest that it is not an equally good choice
 (see~\cite{patera_rozza}). So it is desirable to better understand the
 actual outcome of the RB greedy algorithm used offline.
The concept of {\it greedy} algorithm
appears in many numerical approaches for problems of approximation.
It typically consists in a recursive procedure 
approximating an optimal solution to a complex problem, using a sequence
of  sub-optimal solutions 
incrementally improved. Otherwise stated, each iteration takes the solution of the previous iteration 
as an initial guess and improves it.
In the theory of approximation of functions
in particular~\cite{devore-93,temlyakov-08},
greedy algorithms are used to incrementally compute
the combinations of  functions from a given dictionnary
which best approximate some given function.
The RB greedy algorithm has a somewhat different viewpoint:
it incrementally computes for integers~$N$
some basis functions $u_\N(\mu_n)$, $n=1,\ldots,N$,
spanning a linear space $X_{\N,N}$
that best approximates a  \emph{family} of functions
$u_\N(\mu)$, $\forall \mu\in\Lambda$. The RB greedy algorithm however has a flavour similar to other greedy algorithms
that typically build best-approximants in general classes of
functions.
It is therefore possible to better understand the RB greedy algorithm
using classical ingredients of approximation theory. The notion of
\emph{Kolmogorov width} is an instance of such a classical
ingredient. We refer to~\cite[Chapter 3]{boyaval-these} and~\cite{BMPPT10} for more details and some elements of
analysis of the RB greedy algorithm.

\subsubsection{Extensions of the approach to cases more general
  than~\eqref{eq:PB00}.}

The RB approach of course does not only apply to simple situations
like~\eqref{eq:PB00}. Many more general  situations may be addressed, the major
limitation to the genericity of the approach being the need for constructing  fast computable {\it a posteriori} error estimators.

Instances of problems where the RB approach has been successfully tested
are the following: affine formulations, non-coercive linear
elliptic problems, non-compliant linear elliptic problems, problems with
non-affine
parameters, nonlinear elliptic problems, semi-discretized (nonlinear) parabolic
problems. The purpose of this section is to briefly review these
extensions of our above simple setting. In the next section, we will
then introduce a problem with random coefficients. For simplicity, we
take it almost as simple as the above problem~\eqref{eq:PB00},
see~\eqref{eq:strong}-\eqref{BCbiRB} below. We
anticipate that, if
they involve a random component,  most of the extensions outlined in the present section
could also, \emph{in principle},  be treated using the RB approach.

\paragraph{Affine formulations.} Beyond the simple case presented above in Section~\ref{ssec:outlineRB},
which involves an elliptic operator in
divergence form affinely depending on the parameter,
 the RB approach can be extended to general elliptic problems with variational formulation of the form 
\begin{equation}
\label{eq:general_weak_form}
\text{Find $u(\mu)\in X$ solution to } g(u(\mu),v;\mu)=0\,,\ \forall v\in X\,,
\end{equation}
where the form $g(\cdot,\cdot;\mu)$ on $X\times X$ admits an {\it
  affine} parametrization, that is,  writes
\begin{equation} \label{eq:affineoperator}
g(w,v;\mu) = \sum_{q=1}^Q\Theta_q(\mu)g_q(w,v)\,, \quad \forall w,v \in X\,,\ 
\forall \mu\in\Lambda \,,
\end{equation}
with   parameter-independent forms $\left(g_q(\cdot,\cdot)\right)_{1\le
  q\le Q}$ (where some of the $g_q$ may only depend on $v$) and coefficients
$\left(\Theta_q(\mu)\right)_{1\le q\le Q}$. We emphasize that the whole RB algorithm presented in the simple case above directly translates in this situation. In particular, the matrices used in the online evaluation procedure can be constructed offline.

\paragraph{Non-coercive symmetric linear elliptic problems.} The RB approach can be extended to the case where
the symmetric continuous bilinear form $a(\cdot,\cdot;\mu)$ 
is not coercive but only {\it inf-sup stable}.
An example is the Helmholtz problem treated in~\cite{Pat06:naturalNorms}.
Our discussion of the elliptic problem above can be adapted in a
straightforward way, the only change in offline and online computations 
being that the inf-sup stability constant on $X_\N$:
\begin{equation}
\label{eq:infsup}
0 < \beta_{LB}(\mu) \le \beta(\mu) := 
\underset{w\in X_\N\backslash\{0\}}{\inf}
\underset{v\in X_\N\backslash\{0\}}{\sup}
\frac{a(w,v;\mu)}{\|w\|_X\|v\|_X} \,,\ 
\forall \mu\in\Lambda\,,
\end{equation}
is substituted for $\alpha_{LB}(\mu)$. In practice, the evaluation of $\beta_{LB}(\mu)$ is typically more involved than the evaluation of the coercivity constant $\alpha_{LB}(\mu)$. We refer to~\cite{huyhn-knezevic-chen-hesthaven-patera-09} for an appropriate technique.

\paragraph{Non-compliant linear elliptic problems.}
In~\eqref{eq:outRBput00}, 
the particular choices of $F=l$ for the output
and of symmetric matrices $\uuline{A}(\mu)$ for the definition of the 
 bilinear form $a(\cdot,\cdot;\mu)$
correspond to a particular class of problems called, we recall, {\it compliant}.
Non-compliant linear elliptic problems can be treated as well, but this
is somewhat more technical. These are the cases where, for some $\mu\in\Lambda$ at least,
either $u(\mu)$ is solution to a weak form~\eqref{eq:general_weak_form}
with $g(v,w;\mu)=a(v,w;\mu)-l(w)$, $\forall v,w\in X$
and the bilinear form $a(\cdot,\cdot;\mu)$ is not symmetric,
or the output is $s(\mu)=F(u(\mu))\neq l(u(\mu))$ 
with any linear continuous function $F:X\to\R$.

For instance, we explain how to treat the case of a bilinear form $a(\cdot,\cdot;\mu)$
that is not symmetric, but of course still continuous and inf-sup stable.
The analysis requires considering the solution to the \emph{adjoint} problem
\begin{equation}
\label{eq:linear_dual_weak_form}
\text{Find $\psi(\mu)\in X$ solution to } a(v,\psi(\mu);\mu)=-F(v)\,,\ \forall v\in X\,,
\end{equation}
 along with the corresponding Galerkin discretization $\psi_\N(\mu)\in
X_\N$,  the approximation space $X_{\N,N^\star}^\star$ for the
solution to~\eqref{eq:linear_dual_weak_form}, and an additional RB approximation space $X^\star_{\N,N^\star}\subset X_\N$ 
of dimension $N^\star\ll\N$ . The {\it a posteriori} estimator obtained is similar to~\eqref{eq:bound11},
and writes 
\begin{equation}\label{eq:bound12}
|s_\N(\mu)-s_{\N,N,N^\star}(\mu)|\le 
\Delta^s_{N,N^\star}(\mu) := 
\frac{\|G(\mu)\:u_{\N,N}(\mu)\|_X\:\|G^\star(\mu)\:\psi_{\N,N^\star}(\mu)\|_X}
{\beta_{LB}(\mu)}\,,
\end{equation}
where  $G^\star$ is defined from the adjoint
problem~\eqref{eq:linear_dual_weak_form} similarly to how $G$ is defined
from the original problem. Notice that we again used the inf-sup stability condition~\eqref{eq:infsup}, which indeed holds true after permutation of the arguments $v$ and $w$ (since we work in a finite dimensional space~$X_\N$), the value of the inf sup constant however being not the same. To build the reduced basis of the primal (respectively the dual) problem, in the offline stage, the {\em a posteriori} estimator is based on $\|G(\mu)\:u_{\N,N}(\mu)\|_X$ (respectively $\|G^\star(\mu)\:\psi_{\N,N^\star}(\mu)\|_X$). Apart from the above introduction and use of the adjoint problem, the treatment of the non-compliant case then
basically follows the same lines as that of the compliant case.

Notice that a simple, but less sharp, estimate of the error (namely the left-hand side of~\eqref{eq:bound12}) can be obtained as $\displaystyle\left(\sup_{x \in \X_\N} \frac{|F(x)|}{\|x\|_{X}} \right) \,\|G(\mu)\:u_{\N,N}(\mu)\|_X$. This simple error bound does not involve the solution of any dual problem, and may be of interest in particular in the case when multiple outputs are considered. However, the primal--dual error bound~\eqref{eq:bound12} will be much smaller (since it is quadratic and not linear in the residual) and in many situations, very easy to obtain, since the dual problem is typically simpler to solve than the primal problem (it is indeed linear).

\paragraph{Non-affine parameters.}

We have exploited in several places the affine dependence of $\uuline{A}(\mu)$ in~\eqref{eq:PB00} in terms of the coefficient $\mu$. However, there are many cases for which the dependency on the parameter is more complicated, as for example, when associated with certain kinds of geometric variations.
 Extending the RB approach to the case of non-affine parametrization is
feasible using suitable affine approximations. The computation of  approximations 
$\sum_{m=1}^{M}\beta_m^M(\mu) \uuline{A_m}(x)$
for  functions $\uuline{A}(x;\mu)$ (having in mind as an example the prototypical problem~\eqref{eq:PB00}), is a general problem of approximation.
A  possibility, introduced and further developed in
\cite{barrault04:_empir_inter_method,m2an_magic,maday-nguyen-patera-pau-07} is to modify the standard greedy
procedure described above, using interpolation. In
short, the approach consists in  selecting  the coefficients
$\left(\tilde\beta_m^M(\mu)\right)_{m=1,\ldots,M}$  of  the
approximation 
$
\mathcal{I}_M[g(\cdot;\mu)]:=\sum_{m=1}^{M}\tilde\beta_m^M(\mu)g(\cdot;\mu_m^g)
$ of order $M$ to $g(\cdot;\mu)$ ($g$ denoting here a general bilinear form, as in~\eqref{eq:general_weak_form}--\eqref{eq:affineoperator})  using an interpolation at the so-called  {\it magic points} $x_m$
selected sequentially with 
$$
x_1 \in \mathop{\rm argmax}_{x\in\D} |g(\cdot;\mu_1^g)|\,, \qquad
x_m \in \mathop{\rm argmax}_{x\in\D} 
|g(\cdot;\mu_m^g)-\mathcal{I}_{m-1}[g(\cdot;\mu_m^g)]|\,,
$$
for all $m=2,\ldots,M$. We refer to the contributions cited above for
more details. 

\paragraph{Nonlinear elliptic problems}
For the extension of the RB approach to \emph{nonlinear} problems, one
major difficulty is again the construction of appropriate  {\it a
  posteriori} error estimators, which, additionally, need to be computed
efficiently. Several examples of successful extensions are reported on
in the
literature~\cite{veroy_lions,veroy03:_poster_error_bound_reduc_basis,Patera_Huynh06,deparis07,Patera_Ronquist07,nguyen-09,Calcolo}.
But no general theory can of course be developed in the nonlinear context.

\paragraph{Semi-discretized parabolic problems}

After time-discretization, parametrized
parabolic problems can be viewed as a collection of elliptic problems
with the time variable as an additional parameter. A natural idea is then to build a reduced basis spanned by solutions for given values of the parameter and the time variable. Examples of contributions are~\cite{grepl04:_reduc_basis_approx_time_depen,m2an_magic}. This first approach has been improved by techniques combining the RB idea for the parameter with a proper orthogonal decomposition (POD) in the time variable, first introduced in~\cite{Haasdonk-08} and further discussed in~\cite{knezevic-patera-09}. A route which would be interesting to follow could be to try to adapt on-the-fly, as time goes, the reduced basis which is the most adapted to the current time.

\section{RB Approach for Boundary Value Problems with Stochastic Coefficients}
\label{sec:SPDE}

The first application of the RB approach to a problem with stochastic
coefficients is introduced in~\cite{boyaval-lebris-maday-nguyen-patera-09}. The purpose of this section
is to overview this contribution, in particular showing how the general
RB approach needs to be adapted to the specificities of the problem. We
refer to~\cite{boyaval-lebris-maday-nguyen-patera-09} for all the details omitted below. 

\subsection{Setting  of the problem}

\begin{figure}
\centering
  \includegraphics[scale=0.3]{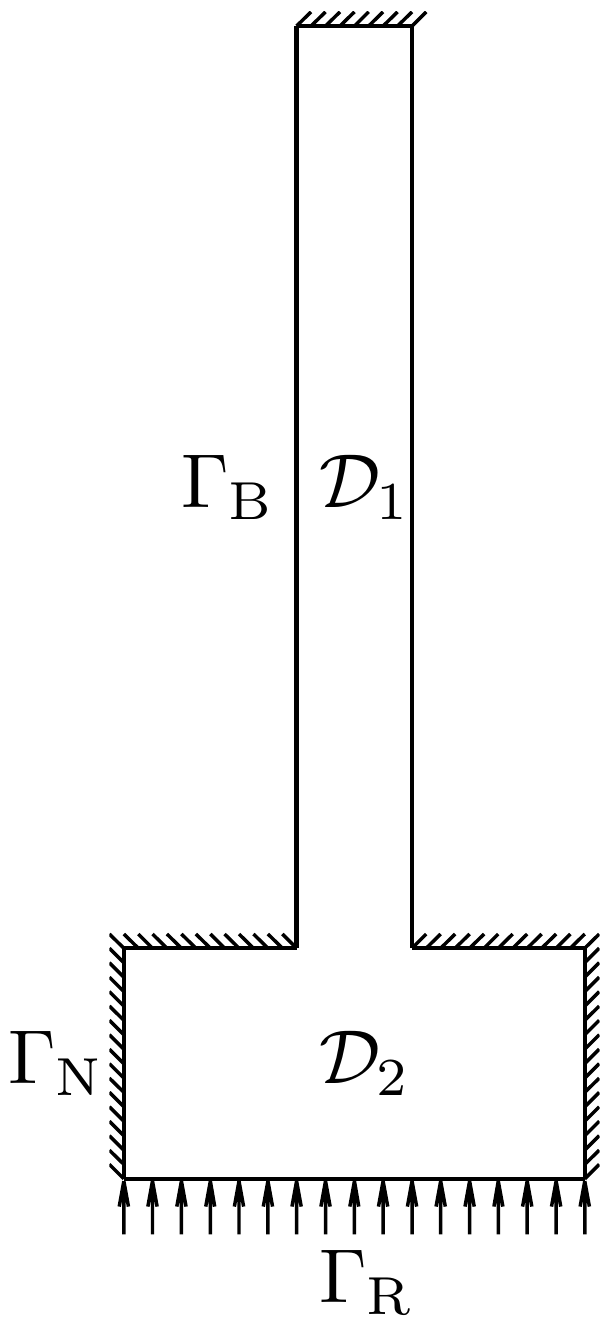}
  \caption{$\D$ has the geometry of a (piece of) heat sink:
a spreader $\mathcal{D}_2$ with a fin $\mathcal{D}_1$ on top.}
  \label{fig0RB}
\end{figure}

Let us denote by $(\Omega,\F,\P)$ a % complete
probability space, and by $\omega\in\Omega$ the stochastic variable.
We consider the stochastic field $U(\cdot,\omega)$ that is  the almost sure  solution to
\begin{equation}
\label{eq:strong}
- {\rm div} \left( \uuline{A}(x) \nabla U(x,\omega) \right) = 0  \,,\ \forall x \in\D\,,
\end{equation}
supplied with a random Robin boundary condition
\begin{equation}
\label{BCbiRB}
\uline{n}(x) \cdot \uuline{A}(x) \nabla U(x,\omega) + B(x,\omega) \: U(x,\omega) = g(x) 
\,,\ \forall x \in\partial\D \ .
\end{equation}
In~\eqref{BCbiRB}, the matrix $\uuline{A}(x)$ writes
$\uuline{A}(x)=\sigma(x)\uuline{I_d}$  where $0<\sigma(x) <\infty$ 
for a.e. $x\in\D$. Of course, $\uline{n}$ denotes  the outward unit normal at the boundary of 
the smooth domain $\D$. The boundary is divided into three non-overlapping open subsets:
$\partial\D= \left(\overline{\Gamma_{\rm N}} \cup
  \overline{\Gamma_{\rm R}} \cup \overline{\Gamma_{\rm B}}\right)$ (see
Fig.~\ref{fig0RB}). The boundary source term $g$ is assumed to vanish
everywhere except on $\Gamma_{\rm R}$ where it has constant unit value:
$g(x)= 1_{\Gamma_{\rm R}}, \forall \: x \in \partial\D$. The scalar random field $B(\cdot,\omega)$, parametrizing the boundary
condition, also vanishes almost everywhere on the boundary $\partial\mathcal{D}$ ,
except on  some subset $\Gamma_{\rm B}$ of the
boundary $\partial\mathcal{D}$ with non-zero measure, where
$0<\bar{b}_{\min}\le B(\cdot,\omega) \le \bar{b}_{\max} <\infty$ almost
surely and almost everywhere. Note that on $\Gamma_{\rm N}$,~\eqref{BCbiRB} thus reduces to homogeneous Neumann conditions.
Physically, 
$U(\cdot,\omega)$ models the steady-state temperature field 
in a heat sink consisting of an isotropic material of thermal
conductivity $\sigma$, contained in the domain $\D$. The sink is
subject to zero heat flux on $\Gamma_{\rm N}$, a
constant flux on $\Gamma_{\rm R}$ modeling the heat source,
and a convective heat transfer on $\Gamma_{\rm B}$.
The Biot number $B$ models the effect of the exterior fluid convection
on the solid thermal conduction problem inside $\D$. In  real world
engineering applications, the value of $B$  is only approximately
known. It is therefore legitimate to encode the uncertainties on $B$ using a  random field $B(\cdot,\omega)$, see~\cite{lienhard}
for more details.

Correspondingly, the solution to \eqref{eq:strong}-\eqref{BCbiRB}, along with any
output computed from this solution, are also random quantities. Only
statitistics on these quantities are relevant. 
We thus consider two \emph{statistical} outputs for the problem:  
the expected value ${\mathbf E}(S)$ and the variance $\mathbf{Var}(S)$
of the random variable 
\begin{equation}
\label{functionalRB}
S(\omega)=\mathcal{F}\left(U(\: \cdot \: ,\omega)\right) = \int_{\Gamma_{\rm R}} U(\: \cdot \: ,\omega) \,
\end{equation}
linearly depending on  the trace of the solution $U(\: \cdot \: ,\omega)$
on $\Gamma_{\rm R}$.

A typical question, example of an  \emph{Uncertainty Quantification}
problem, is to quantify
the sensitivity of the output $S(\omega)$.  Many existing contributions already addressed the issue:
\cite{babuska-tempone-zouraris-05,Deb01,matthies-keese-04,Debusschere04}. 

A possible approach (which we will indeed adopt here) is to evaluate ${\mathbf E}(S)$ and $\mathbf{Var}(S)$
with the plain Monte-Carlo  method
using $M$ independent random variables $(S^m)_{1\leq m\leq M}$
with the same distribution law as $S$. The expectation and the variance
are respectively approached by the empirical sums
\begin{equation}
\label{eq:empirical-estimationsRB}
E_M[(S^m)] = \frac{1}M \sum_{m=1}^M S^m, \qquad
V_M[(S^m)] = \frac{1}{M-1} \sum_{n=1}^M \left(S^n - E_M[(S^m)]\right)^2 \ ,
\end{equation}
where the normalization factors used (respectively $\displaystyle
{1\over M}$ and $\displaystyle
{1\over {M-1}}$) allow, as is traditional in the community of Monte-Carlo
methods, to have unbiased estimators: $\displaystyle {\mathbf E}\left(E_M[(S^m)]\right)={\mathbf
  E}(S) $ and $\displaystyle \mathbf{E}\left(V_M[(S^m)]\right)=\mathbf{Var}(S) $ for all $M$.
Large values of $M$ are typically needed to obtain
from~\eqref{eq:empirical-estimationsRB} accurate
approximations of ${\mathbf E}(S)$ and $\mathbf{Var}(S)$. Since, for
each $m=1,\ldots,M$, a new realization of the random parameter $B$ is
considered and the boundary value problem
\eqref{eq:strong}-\eqref{BCbiRB} has to be solved, the task is clearly
computationally demanding. It is a  many-query context, appropriate for
the application of the RB approach.

% well suited for the deployment of a RB approach if we can deal efficiently
% with the non-affine parameter $b$, and
% outputs ${\mathbf E}(S)$ and $\mathbf{Var}(S)$ (in fact, $E_M[S]$ and $V_M[S]$)
% that are sums over many realizations of the parameter.

\subsection{Discretization of the problem}

We begin by considering the Karhunen--Lo\`{e}ve (abbreviated as KL)
expansion
\begin{equation}
\label{eq:genuineKLRB}
B(x,\omega)
= \overline{b}\,G(x) + \overline{b} \,\sum_{k=1}^{\cal K} \Phi_k(x)\: Y_k(\omega)
\end{equation}
of the
coefficient $B(x,\omega)$
(see~\cite{k.46:_zur_spekt_prozes,loeve-78,schwab-todor-06}).
In~\eqref{eq:genuineKLRB}, $\cal K$ denotes the (possibly
infinite) rank 
of the covariance operator for $B(\cdot,\omega)$, which has eigenvectors $(\Phi_k)_{1\le k\le \mathcal{K}}$
and eigenvalues $(\lambda_k)_{1\le k\le \mathcal{K}}$ (sorted in
decreasing order). The random variables $\left(Y_k\right)_{1\le k\le \mathcal{K}}$ 
are mutually uncorrelated in $L^2_\P(\Omega)$ with zero mean, $G$ is supposed to be normalized $\int_{\partial\D} G =1$
and 
$\overline{b}=\int_{\Omega}d\P(\omega)\int_{\partial\D}B(\cdot,\omega)$
is a fixed intensity factor. 

Based on \eqref{eq:genuineKLRB}, we  introduce the deterministic function 
\begin{equation}
\label{eq:truncatedKLRB} 
b (x,{y}) = \overline{b} \, G(x) 
+ \overline{b} \,\sum_{k=1}^{\mathcal{K}} \Phi_k(x) {y}_k
\end{equation}
defined for almost all~$x \in\partial\D$ and all  ${y}\in\Lambda^{{y}}\subset\mathbb{R}^\mathcal{K}$,
 where $\Lambda^{{y}}$ denotes the range of the sequence 
$Y=\left(Y_k\right)_{1\le k\le \mathcal{K}}$
of random variables appearing in~\eqref{eq:genuineKLRB}. Notice that $B(x,\omega)=b(x,y(\omega))$.

It is next useful to consider, for any positive integer $K\le\mathcal{K}$, 
 truncated versions of the expansions above, and to define, with obvious
 notation,  $U_K(\cdot,\omega)$ as the solution to the
 problem~(\ref{eq:strong})-(\ref{BCbiRB}) where  $B(\cdot,\omega)$  
is replaced by the truncated KL expansion
$B_K(\cdot,\omega)$ at order $K$. Similarly, for all
${y}^K\in\Lambda^{{y}}$, ${u}_K(\,\cdot \,;{y}^K)$  is
defined as the solution to 
\begin{equation}
\label{eq:newRB}
 \left\{
\begin{array}{l}
- {\rm div} \left( \uuline{A}(x) \nabla {u}_K(x;{y}^K) \right) = 0 \,,\ \forall x \in \D \ ,
\\[2ex]
\uline{n}(x)\cdot \uuline{A}(x) \nabla {u}_K(x;{y}^K) 
+ b_K (x,{y}^K) {u}_K(x;{y}^K) = g(x) \,,\ \forall x \in \partial\D,
\end{array} \right.
\end{equation}
where $b_K$ is the $K$-truncated sum~\eqref{eq:truncatedKLRB}.

For a given integer $K\le\mathcal{K}$,
we approximate the random variable $S(\omega)$ 
by $S_K(\omega):=\mathcal{F}\left(U_K(\cdot,\omega)\right)$ where
$U_K(\cdot,\omega)\equiv{u}_K(\cdot;{Y}^K(\omega))$, 
and the statistical outputs ${\mathbf E}(S_K)$ and $\mathbf{Var}(S_K)$ by the empirical sums
\begin{equation}
\label{eq:outRB}
E_M[(S_K^m)] = \frac{1}M \sum_{m=1}^M S_K^m, \qquad
V_M[(S_K^m)] = \frac{1}{M-1} \sum_{n=1}^M \left(S_K^n - E_M[(S_K^m)] \right)^2 \ ,
\end{equation}
using $M$ independent realizations of the random vector ${Y}^K$.
In practice, ${u}_K(\:\cdot\:;{Y}^K_m)$ is approached
using, say, a finite element approximation  ${u}_{K,\N}(\:\cdot\:;{Y}^K_m)$ with $\N\gg1$ degrees of freedom.
Repeating the task for $M$ realizations of the $K$-dimensional random vector ${Y}^K$ may be overwhelming, and this
is where the RB approach comes into the picture. We now present the
application of the RB approach to solve problem~\eqref{eq:newRB},
parametrized by $y^K\in\Lambda^y$.

In echo to our presentation of Section~\ref{sec:initiation-RB}, note that  problem~\eqref{eq:newRB}  is \emph{affine} in the input parameter $y^K$ 
thanks to the KL expansion \eqref{eq:truncatedKLRB} of $b$,
which decouples the dependence on $x$ and the other variables. To use the RB approach for this problem, we consider $S$ in~\eqref{eq:outRB} as the output of the problem, the parameter being $y^K$ (this parameter takes the values $Y^{K,m}$, $m \in \{1, \ldots, M \}$ being the realization number of $Y^K$) and, as will become clear below, the offline stage is standard. On the other hand, in the online stage, the {\em a posteriori} estimation is completed to take into account the truncation error in $K$ in~\eqref{eq:truncatedKLRB}.

Before we turn to this, we emphasize that we have performed above an
\emph{approximation} of the coefficient~$b$, since we have truncated its
 KL expansion. The corresponding error should be estimated. In addition,
  the problem~\eqref{eq:newRB} after truncation might be ill-posed,
  even though the original problem \eqref{eq:strong}-\eqref{BCbiRB}  is
  well posed.  To avoid any corresponding pathological issue,  we
  consider a stochastic coefficient $b$ having 
a KL expansion~\eqref{eq:truncatedKLRB} 
that is positive for any truncation order $K$ (which is a sufficient condition to ensure the well-posedness of~\eqref{eq:strong}-\eqref{BCbiRB}), and which 
converges absolutely a.e.\ in $\partial\D$ when $K\to\mathcal{K}$. For
this purpose, (i) 
we require for $k=1,\ldots,\mathcal{K}$
a uniform bound $\|\Phi_k\|_{L^\infty(\Gamma_{\rm B})} \le \phi$, (ii) we set $Y_k := \Upsilon \sqrt{\lambda_k} Z_k$
with independent random variables $Z_k$
uniformly distributed in the range $(-\sqrt{3},\sqrt{3})$,
$\Upsilon$ being a positive coefficient, and (iii) we also ask $\sum_{k=1}^\mathcal{K}\sqrt{\lambda_k}<\infty$.
Note that, if $\mathcal K = \infty$, condition (iii) imposes a sufficiently fast decay of the
eigenvalues $\lambda_k$ while $k$ increases.
We will see in Section~\ref{sec:applicationRB}
that this fast decay is also important for the practical success of our RB approach. Of course, (i)-(ii)-(iii) are arbitrary conditions that we impose for simplicity. Alternative settings are possible.

\subsection{Reduced-Basis ingredients}

We know from Section~\ref{sec:initiation-RB} that
two essential ingredients in the RB method are an {\it a posteriori}
estimator and a greedy selection procedure.
Like in most applications of the RB method,
both ingredients have to be adapted to the specificities of the present context.

As mentioned above, the statistical outputs~\eqref{eq:outRB} require
new {\it a posteriori} estimators.
Moreover,
the statistical outputs can only be computed after $M$ queries ${Y}^K_m$,
$m=1,\ldots,M$, in the parameter $y^K$,
so  these new {\it a posteriori} estimators 
cannot be used in the offline step.

The global error consists of two, independent contributions: the first
one is related to the  RB
approximation, the second one is related to  the truncation of the KL expansion. 

In the greedy algorithm, we use a standard {\it a posteriori} estimation 
$|S_{K,\N}^m-S_{K,\N,N}^m|\le\Delta^s_{N,K}(Y^K_m)$
for the error between the finite element approximation 
$S_{K,\N}^m:=\mathcal{F}(u_{K,\N}(\:\cdot\:;Y^K_m))$ 
and the RB approximation
$S_{K,\N,N}^m:=\mathcal{F}(u_{K,\N,N}(\:\cdot\:;Y^K_m))$ 
of $S_K^m$ % := \mathcal{F}\left({u}_K(\:\cdot\:;{Y}^K_m)\right)$
at a fixed truncation order $K$,
for any realization $Y_m^K\in\Lambda^y$. %, $m=1,\ldots,M$.
This is classical~\cite{boyaval-08,nguyen-veroy-patera-05,Rozza08:arcme}
and similar to our example of Section~\ref{sec:initiation-RB},
see~\cite{boyaval-lebris-maday-nguyen-patera-09} for details.
Note however that the coercivity constant of the bilinear form for the  variational formulation
\begin{align}
&\text{Find $u(\cdot;y^K)\in H^1(\D)$ s.t. }\nonumber\\
&\int_{\D} \sigma \nabla u(\cdot;y^K) \cdot \nabla v
+ \int_{\Gamma_{\mathrm{B}}} 
b(\cdot,y^K) u(\cdot;y^K) v
=
\int_{\Gamma_{\mathrm{R}}} g\, v, 
\forall \: v \in H^1(\D)\:\label{eqvarform}
\end{align}
of problem~\eqref{eq:newRB} depends on $K$. 
To avoid the additional computation of the coercivity constant for each $K$,
we impose 
$b (x,y^K)\ge\overline{b}G(x)/2$, for all $x\in\Gamma_B$, and thus get a uniform lower bound for the coercivity constant.
In practice, this imposes a limit $0<\Upsilon\le\Upsilon_{\rm max}$
on the intensity factor in the ranges of the random variables $Y_k$,
thus on the random fluctuations of the stochastic coefficient, where $\Upsilon_{\rm max}$ is fixed for all $K \in \{0, \ldots, {\mathcal K}\}$

Let us now discuss the online {\em a posteriori} error estimation. As for the truncation error, an {\it a posteriori} estimation 
$|S_{\N}^m-S_{K,\N}^m|=
|\mathcal{F}(u_{\N}(\:\cdot\:;Y_m^K))-\mathcal{F}(u_{K,\N}(\:\cdot\:;Y_m^K))| 
\le \Delta_{N,K}^{t}(Y_m^K)$
is  derived in~\cite{boyaval-lebris-maday-nguyen-patera-09}.
The error estimators  $\Delta^s_{N,K}(Y^K_m)$ and
$\Delta_{N,K}^{t}(Y^K_m)$, respectively for the RB approximation and the
truncation, 
are eventually combined for $m=1,\ldots,M$
to yield global error bounds  
in the Monte-Carlo estimations of the statistical outputs:
$|E_M[(S_{K,\N,N}^m)]-E_M[(S_{\N}^m)]|\le\Delta_E((S_{K,\N,N}^m))$ 
and $|V_M[(S_{K,\N,N}^m)]-V_M[(S_{\N}^m)]|\le\Delta_V((S_{K,\N,N}^m))$.
The control of the truncation error may be used to improve the performance of the reduced basis method. In particular, if the truncation error happens to be too small compared to the RB approximation error, the truncation rank ${\mathcal K}$ may be reduced.

\subsection{Numerical results}
\label{sec:applicationRB}

Our numerical simulations presented
in~\cite{boyaval-lebris-maday-nguyen-patera-09} are performed on the steady heat conduction problem~(\ref{eq:strong})-(\ref{BCbiRB})
inside the T-shaped heat sink $\D\subset\overline{\D_1}\cup\overline{\D_2}$ 
pictured in Fig.~\ref{fig0RB}.
The heat sink comprises a $2 \times 1$ rectangular substrate (spreader) 
$\D_2 \equiv (-1,1) \times (0,1)$
and a $0.5 \times 4$ thermal fin $\D_1 \equiv (-0.25,0.25) \times (1,5)$ on top.
The diffusion coefficient is piecewise constant,
$\sigma=1_{\D_1} + \sigma_0\:1_{\D_2}$,
where $1_{\D_i}$ of course denotes  the characteristic function of domain $\D_i$ ($i=1,2$).
% We show in Fig.~\ref{fig:spde} 
% numerical results obtained with a regular
The finite element approximation is computed using quadratic finite
elements on a regular mesh, with  $\N=6\,882$ degrees of freedom. The
thermal coefficient is $\sigma_0=2.0$.
To construct the random input field $b(\cdot,\omega)$, we consider the covariance function $\mathbf{Covar}(b(x,\omega)b(y,\omega))=(\overline{b}\Upsilon)^2
\exp(-(x-y)^2/\delta^2)$ for  $\overline{b} = 0.5$,  $\Upsilon=0.058$,
and a  correlation length $\delta=0.5$. We
perform its KL expansion and keep only the largest ${\mathcal K}=25$
terms. We then fix $G(x)\equiv1$ and the variables $Y_k(\omega)$, $1\leq
k\leq{\mathcal K}$, as independent, uniformly distributed random variables. This
defines $b(\cdot,\omega)$ as the right-hand side of \eqref{eq:genuineKLRB}. 

\begin{figure}
\centering
   \includegraphics[trim = 10mm 0mm 0mm 0mm, clip, scale=.4]{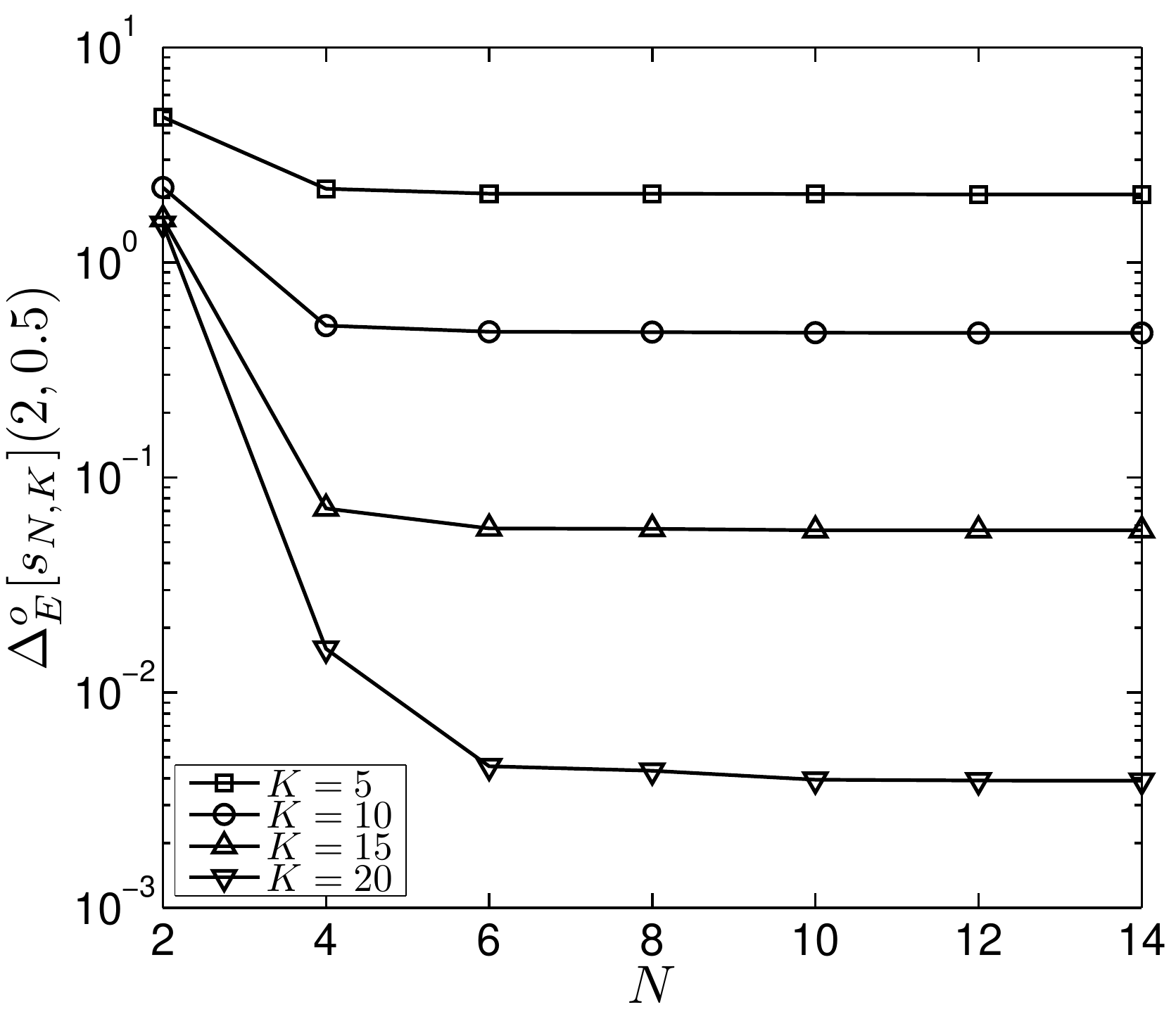}
   \includegraphics[trim = 10mm 0mm 0mm 0mm, clip, scale=.4]{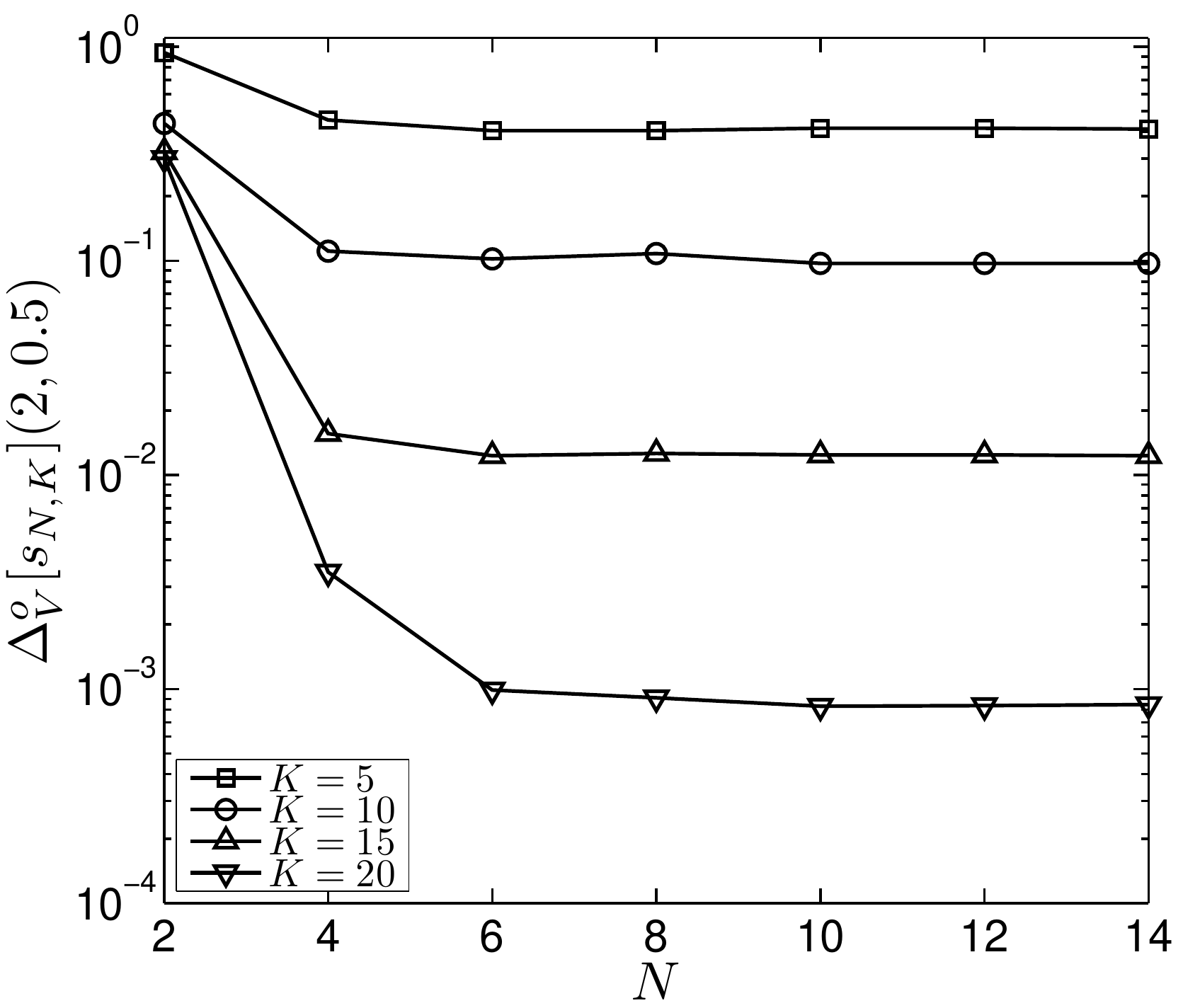}
\caption{Global error bounds for the RB approximation error
and the KL truncation error of
the output expectation (top: $\Delta_E((S_{K,\N,N}^m))$)
and of the output variance (bottom: $\Delta_V((S_{K,\N,N}^m))$),
as functions of the size $N=2,\ldots,14$ of the reduced basis,
at different truncation orders $K=5,10,15,20$.
}\label{fig:spde}
\end{figure}

After computing the reduced basis offline 
with our RB greedy algorithm on a trial sample of size
$|\Lambda_{\rm trial}|=10\,000$,
the global approximation error in the output Monte-Carlo sums 
$E_M[(S_{K,\N,N}^m)]$ and $V_M[(S_{K,\N,N}^m)]$
decays very fast (in fact, exponentially)
with the size $N=1,\ldots,14$ of the reduced basis, 
see Fig.~\ref{fig:spde} with $K=20$.
Note that $M=10\,000$ for the Monte-Carlo sums. We would also like to mention that these reduced bases have actually been obtained letting varying not only the parameter $Y^K$, but also additional parameters (namely the diffusion coefficient $\sigma$ and the mean $\overline{b}$ of the Biot number) but this does not influence qualitatively the results presented here, and we omit this technical issue for simplicity (see~\cite{boyaval-lebris-maday-nguyen-patera-09} for more details).

It is observed that 
the global approximation error for truncated problems at a fixed order $K$ and for various $N$ (the size of the reduced basis)
is quickly dominated by the truncation error. More precisely, 
beyond a critical value $N \ge N_{\rm crit}(K)$, where $N_{\rm crit}(K)$ is increasing with $K$, the global approximation error becomes constant.
Notice that the approximation error is estimated as
usual by {\it a posteriori} estimation techniques.

When $\mathcal{K}$ is infinite (or finite but huge),
the control of the KL truncation error may be difficult. This is a
general issue for problems involving a  decomposition of the stochastic
coefficient. 
Our RB approach is still efficient in some regimes with large $K$, but
not all. 
In particular,
a fast decay of the ranges of the parameters $(y_k)_{1\le k\le K}$
facilitates the exploration of $\Lambda^y$ by the greedy algorithm,
which allows in return to treat large $K$
when the eigenvalues $\lambda_k$ decay sufficiently fast with $k$.

In~\cite{boyaval-lebris-maday-nguyen-patera-09},
we have decreased the correlation length to $\delta=0.2$
and could treat up to $K=45$ parameters, obtaining the results in a
total computational time still fifty times as short as for direct finite element computations.

\section{Variance Reduction using an RB approach}
\label{sec:SDE}

In this section, we present a variance reduction technique based upon an RB approach, which has been proposed recently in~\cite{boyaval-lelievre-09}.
In short, the RB approximation is used as a control variate to reduce the variance of the original Monte-Carlo calculations. 

\subsection{Setting of the problem}

Suppose we need to compute repeatedly, for many values of the parameter $\lambda \in \Lambda$, the Monte-Carlo approximation (using an empirical mean) of the expectation ${\mathbf E}({Z^{\lambda}})$ of a functional 
\begin{equation}\label{outputfunctionalRB}
Z^{\lambda}=g^\lambda(X_T^{\lambda})-\int_0^T f^\lambda(s,X_s^{\lambda})\, ds
\end{equation}
of the solutions $\left( X_t^{\lambda} , t\in[0,T] \right)$ 
to the Stochastic Differential Equation (SDE)
\begin{equation}\label{eq:pbRB}
X_t^{\lambda} = x + \int_0^t b^\lambda(s,X_s^{\lambda})\, ds 
	+ \int_0^t  \sigma^{\lambda}(s,X_s^{\lambda}) dB_s,
\end{equation}
where  $\left( B_t \in \R^d , t\in[0,T] \right)$ is a $d$-dimensional standard Brownian motion.
The parameter $\lambda$ parametrizes the functions $g^\lambda,\, f^\lambda,\, b^\lambda$ and $\sigma^\lambda$.
In~\eqref{eq:pbRB}, we assume $b^\lambda$ and $\sigma^\lambda$ allow for the  It\^o processes $\left(X_t^{\lambda}\in \R^d,t\in[0,T]\right)$ to be well defined, for every $\lambda \in \Lambda$.
Notice that we have supplied the equation with the deterministic initial condition $X^\lambda_0 = x\in \R^d$.
In addition, $f^\lambda$ and $g^\lambda$ are also assumed smooth,
such that $Z^{\lambda}\in L^2(\Omega)$. Recall that a symbolic concise notation for~\eqref{eq:pbRB} is
$$dX_t^{\lambda} = b^\lambda(t,X_t^{\lambda})\, dt 
	+  \sigma^{\lambda}(t,X_t^{\lambda}) \, dB_t \text{ with } X_0^\lambda=x.$$

Such parametrized problems are encountered in numerous applications,
such as the calibration of the volatility in finance,
or the molecular simulation of Brownian particles in materials science.
For the applications in finance, ${\mathbf E}(Z^\lambda)$ is typically the price of an European option in the Black-Scholes model, and $\lambda$ enters the diffusion term (the latter being called the volatility in this context). The calibration of the volatility consists in optimizing $\lambda$ so that the prices observed on the market are close to the prices predicted by the model. Any optimization procedure requires the evaluation of ${\mathbf E}(Z^\lambda)$ for many values of $\lambda$. On the other hand, the typical application we have in mind in materials science is related to polymeric fluids modelling. There, ${\mathbf E}(Z^\lambda)$ is a stress tensor which enters the classical momentum conservation equation on velocity and pressure, and $X^\lambda_t$ is a vector describing the configuration of the polymer chain, which evolves according to an overdamped Langevin equation, namely a stochastic differential equation such as~\eqref{eq:pbRB}. In this context, $\lambda$ is typically the gradient of the velocity field surrounding the polymer chain at a given point in the fluid domain. The parameter $\lambda$ enters the drift coefficient $b^\lambda$. The computation of the stress tensor has to be performed for each time step, and for many points in the fluid domain, which again defined a many-query context, well adapted to the RB approach. For more details on these two applications, we refer to~\cite{boyaval-lelievre-09}.

\medskip
We consider the general form~(\ref{outputfunctionalRB})--(\ref{eq:pbRB}) 
of the problem and as output the Monte-Carlo estimation ${\rm E_M}[(Z^{\lambda}_m)]=\frac{1}{M}\sum_{m=1}^M Z^{\lambda}_m$ parametrized by $\lambda\in\Lambda$, where we recall $(Z^{\lambda}_m)$ denotes i.i.d. random variables with the same law as $Z^\lambda$. These random variables are build in practice by considering a collection of realizations of~\eqref{eq:pbRB}, each one driven by a Brownian motion independent from the others. In view of the  Central Limit Theorem, the rate at which the Monte-Carlo approximation ${\rm E_M}[(Z^{\lambda}_m)]$ approaches its limit ${\mathbf E}(Z^\lambda)$ is given by $\frac{1}{\sqrt{M}}$, the prefactor being proportional to the variance of $Z^\lambda$. A standard approach for reducing the amount of computations is therefore
\emph{variance
reduction}~\cite{arouna-03,melchior-ottinger-95,ottinger-vandenbrule-hulsen-97,bonvin-picasso-99,hammersley-handscomb-64,milstein-tretyakov-06}.
We focus on one particular variance reduction technique: the control variate method. It consists in introducing a so called {\em control variate} $Y^{\lambda}\in L^2(\Omega)$, assumed centered here for simplicity:
$${\mathbf E}({Y^{\lambda}})=0,$$ and in considering the equality:
$$ {\mathbf E}({Z^{\lambda}}) = {\mathbf E}({Z^{\lambda}-Y^{\lambda}}).$$
The expectation ${\mathbf E}({Z^{\lambda}-Y^{\lambda}})$ is approximated by Monte-Carlo estimations ${\rm E_M}[(Z^{\lambda}_m-Y^{\lambda}_m)]$ which hopefully have, for a well chosen $Y^{\lambda}$, a smaller statistical error than direct Monte-Carlo estimations 
${\rm E_M}[(Z^{\lambda}_m)]$ of ${\mathbf E}({Z^{\lambda}})$. More precisely, $Y^\lambda$ is expected to be chosen so that $\mathbf{Var}(Z^{\lambda})\gg\mathbf{Var}(Z^{\lambda}-Y^{\lambda})$.
The central limit theorem yields
\begin{equation}\label{eq:slln}
{\rm E}_M[(Z^\lambda_m-Y^\lambda_m)]:= \frac1M \sum_{m=1}^M (Z_m^{\lambda}-Y^\lambda_m) \xrightarrow[M\to\infty]{\P-a.s.}
{\mathbf E}(Z^\lambda-Y^\lambda),
\end{equation}
where the error is controlled by confidence intervals, in turns functions of the variance of the random variable at hand. The empirical variance
\begin{equation} \label{eq:variance-estimationRB}
{\rm Var}_M\left((Z^{\lambda}_m-Y^\lambda_m)\right):=\frac{1}{M-1} \sum_{n=1}^{M}
 \left(Z^{\lambda}_n- Y^\lambda_n-{\rm E}_M((Z^{\lambda}_m-\ Y^\lambda_m))\right)^2  
\end{equation}
which, as $M \to \infty$, converges to ${\rm Var}(Z^\lambda)$, yields a computable error bound. The Central Limit Theorem indeed states that: for all $a > 0$,
\begin{equation}
\label{eq:cltRB}
\P\left( \left| {\mathbf E}(Z^{\lambda}-Y^\lambda) - {\rm E}_M\left((Z^{\lambda}_m- Y^\lambda_m)\right) \right|
\le a\sqrt{\frac{{\rm Var}_M\left((Z^{\lambda}_m- Y^\lambda_m)\right)}{M}} \right)
\xrightarrow[M\to\infty]{} \int_{-a}^a \frac{e^{-x^2/2}}{\sqrt{2\pi}} dx  \,.
\end{equation}
Evaluating the empirical variance~\eqref{eq:variance-estimationRB} is therefore an ingredient in Monte-Carlo computations, similar to what {\em a posteriori} estimates are for a deterministic problem.

Of course, the ideal control variate
is, $\forall\lambda\in\Lambda$:
\begin{equation}\label{eq:bestRB}
Y^{\lambda}=Z^{\lambda}-{\mathbf E}(Z^{\lambda})\,,
\end{equation}
since then, $\Var{Z^{\lambda}-Y^{\lambda}}=0$.
This is however not a practical control variate since ${\mathbf E}(Z^{\lambda})$ itself, the quantity we are trying to evaluate, is necessary to compute~\eqref{eq:bestRB}.
 It\^o calculus shows that the optimal control variate~\eqref{eq:bestRB} also writes:
\begin{equation} \label{eq:ideal_controlRB}
Y^{\lambda}= 
\int_0^T \nabla u^\lambda(s,X_s^{\lambda}) \cdot \sigma^{\lambda}(s,X_s^{\lambda}) dB_s,
\end{equation}
where $u^\lambda(t,y)\in C^1\left([0,T],C^2(\R^d)\right)$ satisfies the backward Kolmogorov equation:
\begin{equation}\label{eq:PDERB}
\left\{
\begin{array}{l} \displaystyle 
\partial_t u^\lambda + b^\lambda\cdot\nabla u^\lambda 
	+ \frac{1}{2} \sigma^\lambda(\sigma^\lambda)^T:\nabla^2 u^\lambda = f^{\lambda} \ ,
\\
u^\lambda(T,\cdot) = g^{\lambda}(\cdot).
\end{array}
\right.
\end{equation}
% More explicitly, $\nabla u^\lambda\equiv\nabla_y u^\lambda(t,y)$
% and $\sigma^\lambda(\sigma^\lambda)^T:\nabla^2 
% u^\lambda\equiv\sum_{i,j,k=1}^d
% \sigma^\lambda_{ik}(t,y)\sigma^\lambda_{jk}(t,y)\partial_{y_i,y_j}^2 u^\lambda(t,y) $.
Even using this reformulation, the choice~\eqref{eq:bestRB} is impractical since solving the partial differential equation~\eqref{eq:PDERB} is at least as difficult as computing ${\mathbf E}(Z^{\lambda})$. We will however explain now that both "impractical" approaches above may give birth to a practical variance reduction method, when they are combined with a RB type approximation.

Loosely speaking, the idea consists in: (i) in the offline stage, compute fine approximations of ${\mathbf E}(Z^{\lambda})$ or respectively $u^\lambda$ for some appropriate values of $\lambda$, in order to obtain fine approximations of the optimal control variate $Y^\lambda$ (at those values) and (ii) in the online stage, for a new parameter~$\lambda$, use as a control variate the best  linear combination of the variables built offline.

% \begin{remark}
% The computation of the expectation of~\eqref{outputfunctionalRB}
% can also be achieved with quadrature formulas in a fully deterministic setting
% using the probability density functional of $(X^\lambda_t)$,
% solution to the Fokker-Planck equation 
% associated with the SDE~\eqref{eq:pbRB}. This goes back to the classical RB setting, see~\cite{knezevic-patera-09}.
% The approach here is different since we consider a stochastic discretization,
% typically more suitable for high-dimensional settings $d\ge 4$.
% \end{remark}

\subsection{Two algorithms for variance reduction by the RB approach}

Using suitable time discretization methods~\cite{kloeden-platen-00}, realizations of the stochastic process~\eqref{eq:pbRB} and the corresponding functional~\eqref{outputfunctionalRB}
can be computed for any $\lambda\in\Lambda$, as precisely as needed. Leaving aside all technicalities related to time discretization, we thus focus on the Monte Carlo discretization.

We construct two algorithms, which can be outlined as follows.

{\bf Algorithm~1} (based on formulation~\eqref{eq:bestRB}):
\begin{itemize}
\item {\em Offline stage}: Build an appropriate set of values $\{\lambda_1, \ldots, \lambda_N \}$ and, concurrently, for each $\lambda \in \{\lambda_1, \ldots, \lambda_N \}$ compute an accurate approximation ${\rm E}_{M_{\rm large}}[(Z^{\lambda}_m)]$ of ${\mathbf E}(Z^{\lambda})$ (for  a very large number $M_{\rm large}$ of realizations). At the end of the offline step, accurate approximations $$\tilde{Y}^{\lambda} = Z^{\lambda}-{\rm E}_{M_{\rm large}}[Z^{\lambda}_m]$$ 
of the optimal control variate  $Y^\lambda$ are at hand. The set of values $\{\lambda_1, \ldots, \lambda_N \}$ is chosen in order to ensure the maximal variance reduction in the forthcoming online computations (see below for more details). 
\item {\em Online stage}: For any $\lambda \in \Lambda$, compute a control variate $\tilde{Y}^\lambda_N$ for the Monte-Carlo estimation of ${\mathbf E}(Z^{\lambda})$ as a linear combination of $$(\bar Y_i=\tilde{Y}^{\lambda_i})_{1 \le i \le N}.$$ 
\end{itemize}

{\bf Algorithm~2} (based on formulation~\eqref{eq:ideal_controlRB}):
\begin{itemize}
\item {\em Offline stage}: Build an appropriate set of values $\{\lambda_1, \ldots, \lambda_N \}$ and, concurrently, for each $\lambda \in \{\lambda_1, \ldots, \lambda_N \}$, compute an accurate approximation $\tilde{u}^\lambda$ of $u^\lambda$, by solving the partial differential equation~\eqref{eq:PDERB}. The set of values $\{\lambda_1, \ldots, \lambda_N \}$ is chosen in order to ensure the maximal variance reduction in the forthcoming online computations  (see below for more details).
\item {\em Online stage}: For any $\lambda \in \Lambda$, compute a control variate $\tilde{Y}^\lambda_N$ for the Monte-Carlo estimation of ${\mathbf E}(Z^{\lambda})$ as a linear combination of $$\left(\bar Y_i=\int_0^T \nabla \tilde u^{\lambda_i}(s,X_s^{\lambda}) \cdot \sigma^{\lambda}(s,X_s^{\lambda}) dB_s
\right)_{1 \le i \le N}.$$ 
\end{itemize}
In both algorithms, we denote by $\tilde Y^\lambda_N$ the control variate built online as a linear combinations of the $\bar Y_i$'s, the lowerscript index $N$ emphasizing that the approximation is computed on a basis with $N$ elements.
An important practical ingredient in both algorithms is to use for the computation of $Z^\lambda$ the exact same Brownian motions as those used to build the control variates.

The construction of set of values $\lambda \in \{\lambda_1, \ldots, \lambda_N \}$ in the offline stage of both algorithms is done using a greedy algorithm similar to those considered in the preceeding sections. The only difference is that the error estimator used is the empirical variance. Before entering that, we need to make precise how the linear combinations are built online, since this linear combination construction is also used offline to choose the $\lambda_i$'s.

The online stage of both algorithms 1 and 2 follow the same line: for a given parameter value $\lambda \in \Lambda$, a control variate $\tilde Y^\lambda_N$ for $Z_\lambda$ is built as an appropriate linear combination of the control variates $(\bar Y_i)_{1 \le i \le N}$ (obtained from the offline computations). The criterium used to select this appropriate combination is based on a minimization of the variance:
\begin{equation}\label{eq:minimumY}
\tilde Y^\lambda_N=\sum_{n=1}^N \alpha_n^* \, {\bar Y}_n,
\end{equation}
where
\begin{equation}\label{eq:minimumY}
(\alpha_n^*)_{1 \le n \le N}= \arg\min_{(\alpha_n)_{1 \le n \le N} \in \R^N} \mathbf{Var}\left(Z^{\lambda} - \sum_{n=1}^N \alpha_n {\bar Y}_n\right).
\end{equation}
In practice the variance in~\eqref{eq:minimumY} is of course replaced by its empirical approximation ${\rm Var}_{M_{\rm small}}$.   Notice that we have an error estimate of the Monte Carlo approximation by considering ${\rm Var}_{M_{\rm small}}\left(Z^{\lambda}-\tilde{Y}_N^{\lambda}\right)$.
It is easy to check that the least squares problem~\eqref{eq:minimumY} is computationally inexpensive to solve since it amounts to solving a linear $N\times N$ system, with $N$ small. More precisely, this linear system writes:
$$C_{M_{\rm small}} \alpha^* = b_{M_{\rm small}}$$
where $\alpha^*$ here denotes the vector with components $\alpha_n^*$, $C_{M_{\rm small}}$ is a matrix with $(i,j)$-th entry 
$${\rm Cov}_{M_{\rm small}}({\bar Y}_{i,m},{\bar Y}_{j,m})$$
and $b_{M_{\rm small}}$ is a vecteur with $j$-th component
$${\rm Cov}_{M_{\rm small}}(Z^\lambda_{m},{\bar Y}_{j,m})$$
where for two collections of random variables $U_m$ and $V_m$, 
$${\rm Cov}_M(U_{m},V_{m})=\frac{1}{M} \sum_{m=1}^{M} U_{m} V_{m} - \left(\frac{1}{M} \sum_{m=1}^{M} U_m\right) \left(\frac{1}{M} \sum_{m=1}^{M} V_m\right).$$
In summary, the computational complexity of one online evaluation is the sum of the computational cost of the construction of $b_{M_{\rm small}}$ (wich scales like $N M_{\rm small}$), and of the resolution of the linear system (which scales like $N^2$ for Algorithm 1 since the SVD decomposion of $C_{M_{\rm small}}$ may be precomputed offline, and scales like $N^3$ for Algorithm 2, since the whole matrix $C_{M_{\rm small}}$ has to be recomputed for each new value of $\lambda$).

The greedy algorithms used in the offline stages follow the same line as in the classical RB approach. More precisely, for Algorithm 1, the offline stage writes:  Let $\lambda_1 \in \Lambda_{\rm trial}$ be already chosen and compute ${\rm E}_{M_{\rm large}}(Z^{\lambda_1})$. Then, for $i=1, \ldots, N-1$, for all $\lambda \in \Lambda_{\rm trial}$, compute $\tilde{Y}_i^{\lambda}$ and inexpensive approximations:
$$
E_i(\lambda):={\rm E}_{M_{\rm small}}(Z^{\lambda}-\tilde{Y}_i^{\lambda})\text{ for }{\mathbf E}(Z^{\lambda}) \ ,
$$
$$
\epsilon_i(\lambda):={\rm Var}_{M_{\rm small}}\left(Z^{\lambda}-\tilde{Y}_i^{\lambda}\right)\text{ for } 
\Var{Z^{\lambda}-\tilde{Y}_i^{\lambda}} \ .
$$
Select $\lambda_{i+1} \in 
\underset{\lambda \in\Lambda_{\rm trial}\backslash\{\lambda_j,j=1,\ldots, i\}}{\argmax}
\left\{ \epsilon_i(\lambda) \right\} $, and compute ${\rm E}_{M_{\rm large}}(Z^{\lambda_{i+1}})$.

In practice, the number $N$ is determined such that $\epsilon_N(\lambda_{N+1})\le\varepsilon$, for a given threshold $\varepsilon$.
The greedy procedure for Algorithm 2 is similar.

\subsection{Reduced-Basis ingredients}

The algorithms presented above to build a control variate using a reduced basis share many features with the classical RB approach. The approach follows a two-stage offline / online strategy. The reduced basis is built using snapshots (namely solutions for well chosen values of the parameters). An inexpensive error estimator is used both in the offline stage to build the reduced basis in the greedy algorithm, and in the online stage to check that the variance reduction is correct for new values of the parameters. The construction of the linear combinations for the control variates is based on a minimization principle, which is reminiscent of the Galerkin procedure~\eqref{eq:primal00N}.

The practical efficiency observed on specific examples is similar for the two algorithms. They both satisfactorily reduce variance. Compared to the plain Monte Carlo method without variance reduction, the variance is divided at least by a factor $10^2$, and typically by a factor $10^4$. Algorithm~2 appears to be computationally much more demanding than Algorithm~1 and less general, since it requires the computation (and the storage) of an approximation of the solution to the backward Kolmogorov equation~\eqref{eq:PDERB} for a few values of the parameter. In particular, Algorithm~2  seems impractical for high dimensional problems ($X^\lambda_t \in \R^d$ with $d$ large). On the other hand, Algorithm~2 seems to be more robust with respect to the choice of $\Lambda_{\rm trial}$: it yields good variance reduction even for large variations of the parameter~$\lambda$, in the online stage. We refer to~\cite{boyaval-lelievre-09} for more details.

Notice also that Algorithm~1 is not restricted to a random variable $Z^\lambda$ that is defined as a functional of a solution to a SDE. The approach can be generalized to any parametrized random variables, as long as there is a natural method to generate correlated samples for various values of the parameter. A natural setting for such a situation is the computation of a quantity ${\mathbb E}(g^\lambda(X))$ for a random variable $X$ with given arbitrary law, {\em independent of the parameter~$\lambda$}. In such a situation, it is easy to generate correlated samples by using the same realizations of the random variable $X$ for various values of the parameter $\lambda$.

\subsection{Numerical Results}

\begin{figure}[htbp]
\centering
\includegraphics[trim = 8mm 0mm 20mm 10mm, clip, scale=.3]{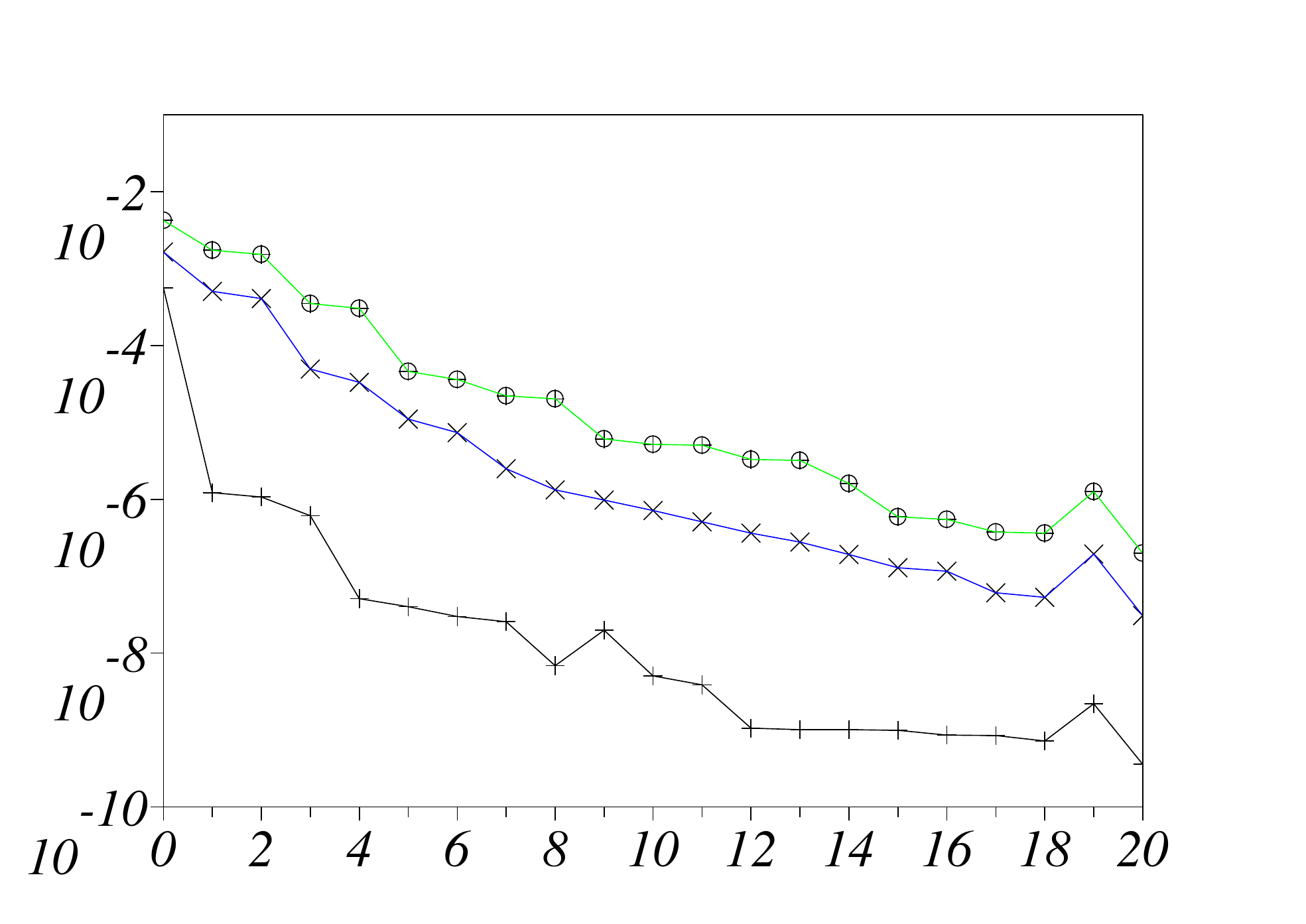}
\includegraphics[trim = 8mm 0mm 20mm 10mm, clip, scale=.3]{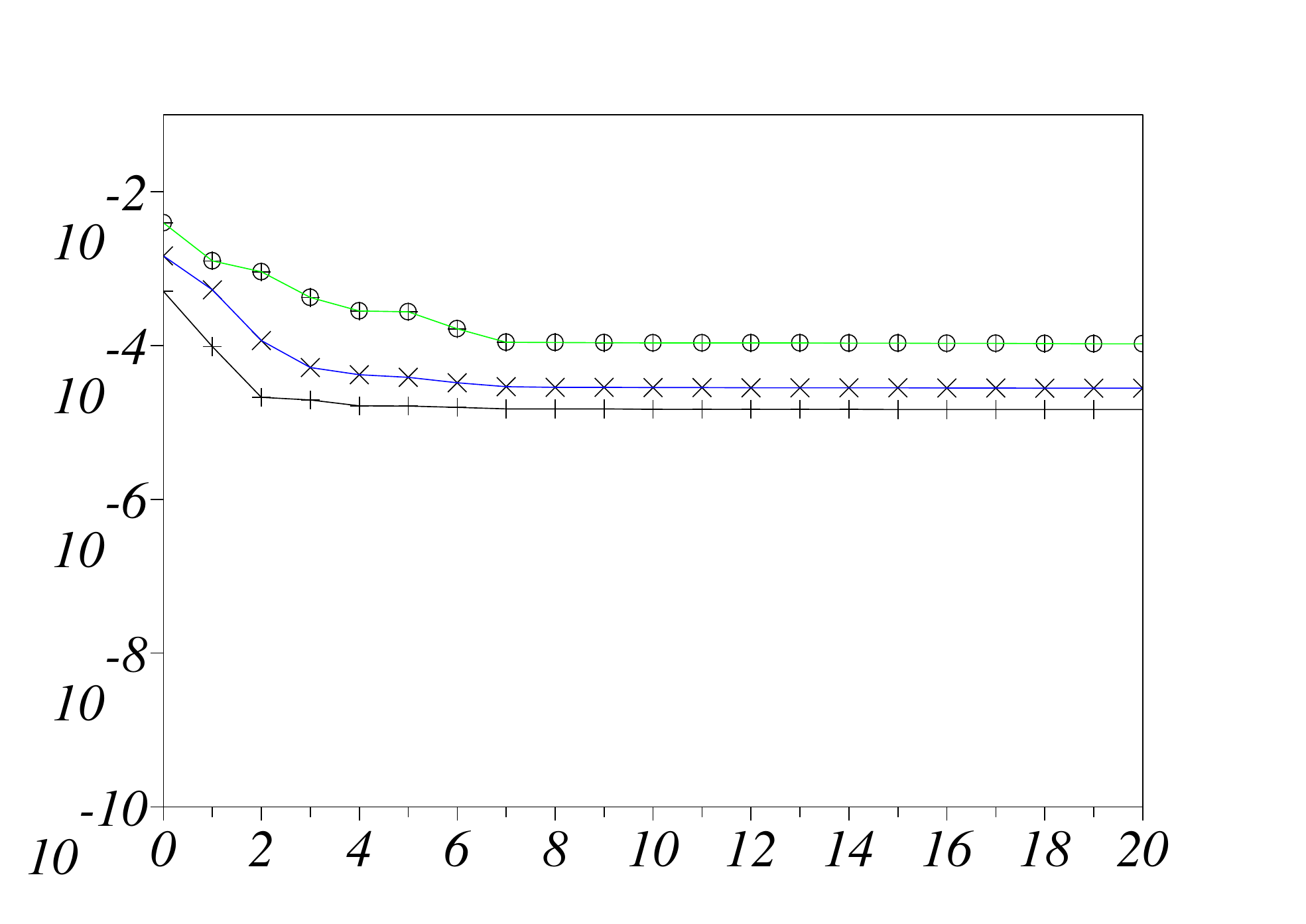}

\includegraphics[trim = 8mm 0mm 20mm 10mm, clip, scale=.3]{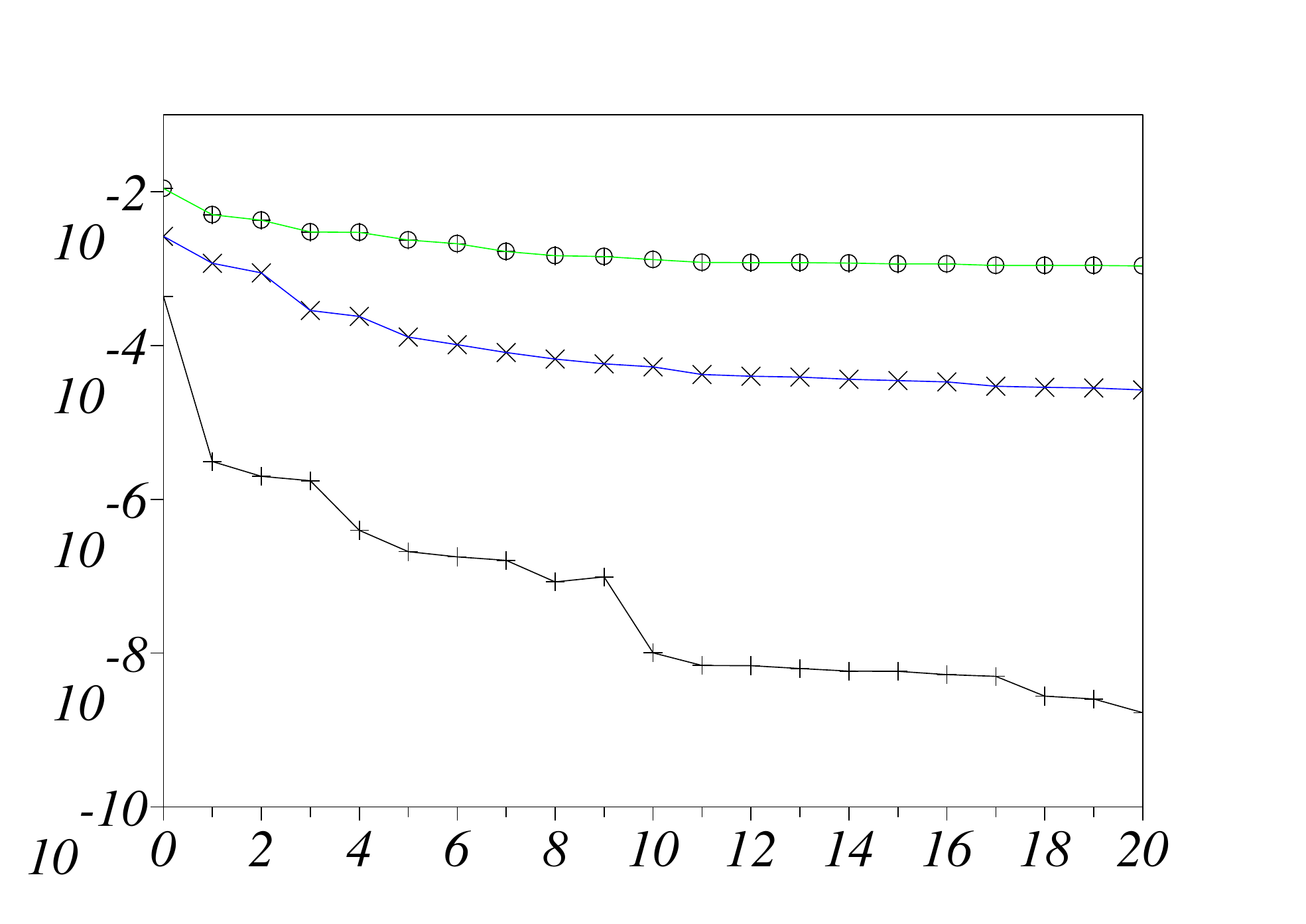}
\includegraphics[trim = 8mm 0mm 20mm 10mm, clip, scale=.3]{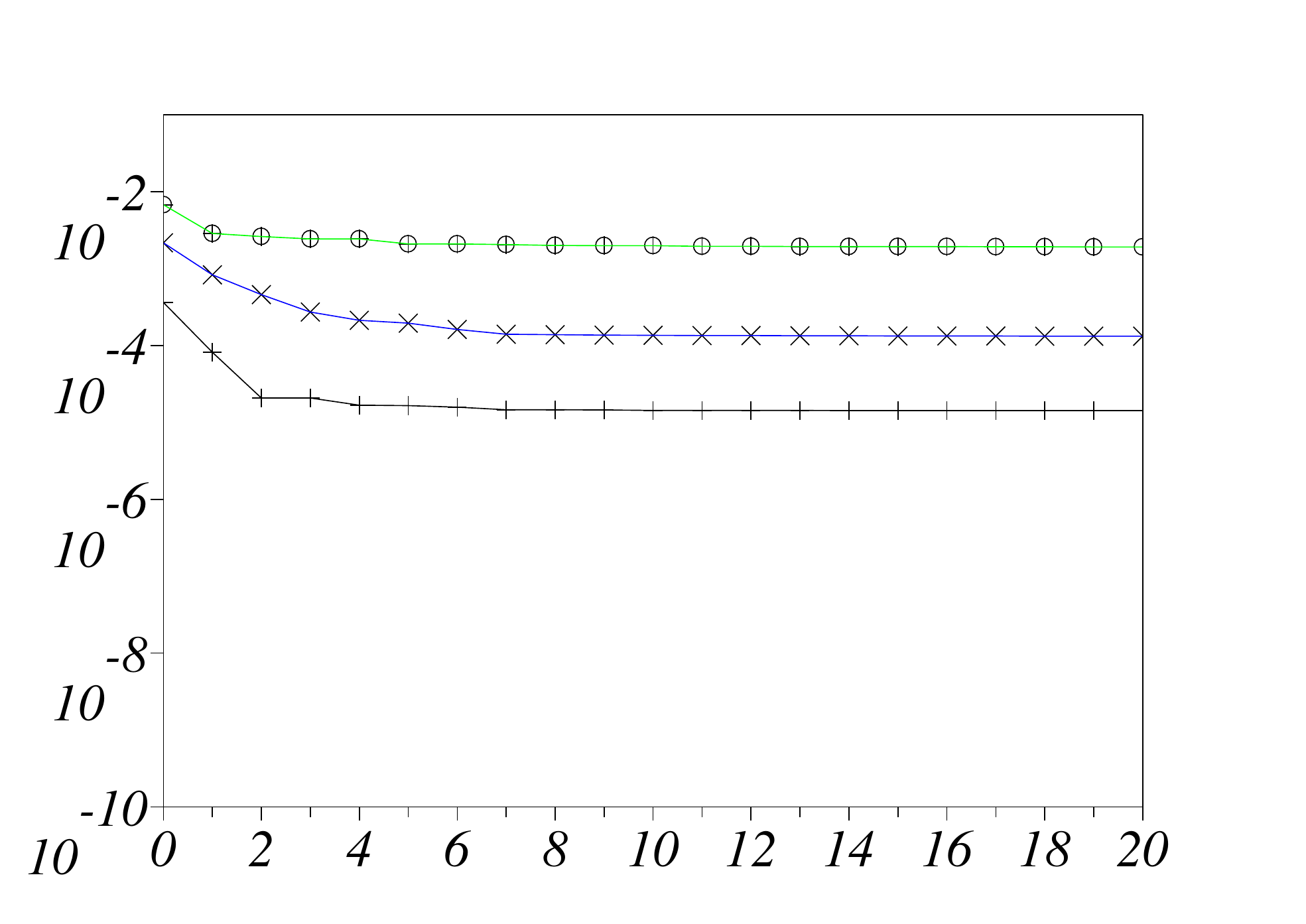}
\caption{Algorithm~1 (left) and 2 (right) for FENE model with $b=16$. The $x$-axis is the size $N$ of the reduced basis.
We represent the minimum $+$, mean $\times$ and maximum $\circ$ of 
${\rm Var}_M[Z^{\lambda} - {\tilde Y}_N^{\lambda}]/{\rm E}_M[Z^{\lambda} - {\tilde Y}_N^{\lambda}]^2$
over online test samples $\Lambda_{\rm test}\subset\Lambda$ (top)
and $\Lambda_{\rm test wide}\supset\Lambda$ (bottom) of parameters.
\label{fig:FENEdistrib6RB}
}
\end{figure}

The numerical results shown on Figure~\ref{fig:FENEdistrib6RB} are taken from~\cite{boyaval-lelievre-09} and relate to the second application mentioned in the introduction, namely multiscale models for polymeric fluids (see~\cite{LL07} for a general introduction). In this context, the non-Newtonian stress tensor is defined by the Kramers formula as an expectation ${\mathbf E}({Z^{\lambda}})$ of the random variable:
\begin{equation}\label{eq:kramersRB}
Z^{\lambda} = {X}_T^{\lambda}\otimes{F}({X}_T^{\lambda}) \,,
\end{equation}
where $X^\lambda_t$ is a vector modelling the conformation of the polymer chain. The latter evolves according to an overdamped Langevin equation:
\begin{equation}\label{eq:langevinRB}
 d{X}_t^{\lambda} 
= \left( {\lambda}
\,{X}^{\lambda}_t-{F}({X}^{\lambda}_t) \right)\, dt + d{B}_t.
\end{equation}
Equation~\eqref{eq:langevinRB} holds at each position of the fluid domain, 
the parameter $\lambda\in\R^{d\times d}$  ($d=2$ or $3$) being the local instantaneous value
of the velocity gradient field at the position considered.
The evolution of the "end-to-end vector" ${X}_t^{\lambda}$ is governed by three forces:
a hydrodynamic force ${\lambda}{X}_t^{\lambda}$,
Brownian collisions ${B}_t$ against the solvent molecules, and 
an entropic force ${F}({X}_t^{\lambda})$ specific to the polymer molecule.
Typically, this entropic force reads either ${F}({X}_t^{\lambda})={X}_t^{\lambda}$ (for the Hookean dumbbells), 
or ${F}({X}_t^{\lambda})=\frac{{X}_t^{\lambda}}{1 - |{X}_t^{\lambda}|^2/b}$
(for the Finitely-Extensible Nonlinear Elastic (FENE) dumbells, assuming $|{X}_t^{\lambda}|<\sqrt{b}$).

The numerical simulations of the flow evolution of a polymeric fluid using such a model
typically consist, on many successive time slots $[nT,(n+1)T]$, of two steps:
(i) the computation of~\eqref{eq:langevinRB}, for a given gradient velocity field $\lambda$, at many points of the fluid domain (think of the nodes of a finite element mesh) and 
(ii) the computation of a new velocity gradient field in the fluid domain, for a given value of the non-Newtonian stress tensor, by solving the classical momentum and mass conservation equations, we omit here for brevity. Thus, ${\mathbf E}({Z^{\lambda}})$ has to be computed for many values $\lambda$ corresponding to many spatial positions and many possible velocity fields at each such positions in the fluid domain.

In the numerical simulations of Figure~\ref{fig:FENEdistrib6RB}, the SDE~\eqref{eq:langevinRB} for FENE dumbbells when $d=2$
is discretized with the Euler-Maruyama scheme
using $100$ iterations with a constant time step $\Delta t = 10^{-2}$
starting from a deterministic initial condition $x=(1,1)$. Reflecting boundary conditions are imposed on the boundary of the ball with radius $\sqrt{b}$.
For $b=16$ and $|\Lambda_{\rm trial}|=100$ trial parameter values
randomly chosen in the cubic range $\Lambda=[-1,1]^3$
(the traceless matrix $\uuline{\lambda}$ has 
entries $(\lambda_{11}=-\lambda_{22},\lambda_{12},\lambda_{21})$),
a greedy algorithm is used to incrementally select $N=20$
parameter values after solving $|\Lambda_{\rm trial}|=100$
least-squares problems~\eqref{eq:minimumY} (with $M_{\rm small}=1000$) 
at each step of the greedy algorithm
(one for each of the trial parameter values $\lambda\in\Lambda_{\rm trial}$).
Then, the $N=20$ selected parameter values are used online 
for variance reduction of a test sample of
$|\Lambda_{\rm test}|=1000$ random parameter values.

The variance reduction obtained online by Algorithm~1 
with $M_{\rm large}=100\:M_{\rm small}$ is very interesting,
of about $4$ orders of magnitude.
For the Algorithm~2, we use the {\em exact} solution $\tilde u^\lambda$ 
to the Kolmogorov backward equation for Hookean dumbells 
as an approximation to $u^\lambda$ solution to~\eqref{eq:PDERB}.
This also yields satisfying variance reduction 
though apparently not as good as in Algorithm~1.
As mentioned above, Algorithm~2 is computationally more demanding but seems to be slightly more robust than Algorithm~1 (namely 
when some online sample test
$\Lambda_{\rm test wide}$ uniformly distributed in $[-2,2]^3$
extrapolates the trial sample used offline, see Fig.~\ref{fig:FENEdistrib6RB}).

% \begin{figure}[htbp]
% \centering
% \includegraphics[trim = 8mm 0mm 20mm 10mm, clip, scale=.37]{minmedmax15r2off.pdf}
% \includegraphics[trim = 8mm 0mm 20mm 10mm, clip, scale=.37]{minmedmax15r1on.pdf}
% \caption{Algorithm~1 for FENE model with $b=9$: 
% Minimum $+$, mean $\times$ and maximum $\circ$ of ${\rm Var}_M[Z^{\lambda} - {\tilde Y}_N^{\lambda}]$
% in samples of parameters
% (left: offline sample $\Lambda_{\rm trial}\setminus\{\lambda_n,n=1,\ldots,N\}$; right: online sample $\Lambda_{\rm test}$) 
% with respect to the size $N$ of the reduced basis.
% \label{fig:FENEdistrib1}
% }
% \end{figure}

Our numerical tests, although preliminary, already show that the reiterated computations of parametrized Monte-Carlo estimations
seem to be a promising opportunity of applications for RB approaches. More generally, even if RB approaches may not be accurate enough for some applications, they may be seen as good methods to obtain first estimates, which can then be used to construct more refined approximations (using variance reduction as mentioned here, or maybe preconditionning based on the coarse-grained RB model). This is perhaps the most important conclusion of the work described in this section.

% Much remains to be done theoretically:
% the situation is similar to the standard RB method.
% And many improvements of the initial work~\cite{boyaval-lelievre-09}
% are also certainly possible on the numerical side,
% possibly dedicated to specific applications.
% Among the possible directions for further exploration of the method,
% high-dimensional problems seem particularly interesting ($d\ge 4$),
% which may require probabilistic approximations of $u^\lambda$
% (see~\cite{newton-94}).

\section{Perspectives}\label{sec:futur}

The standard RB method has proved numerically
efficient and reliable at reducing the cost of computations 
for the approximation of solutions to parametrized boundary value problems
in numerous benchmark many-query frameworks.
These accomplishments claim
for a wider use of the RB method in more realistic settings,
and even suggest that some RB ideas could still be extended
in numerous many-query frameworks yet largely unexplored,
including the stochastic context.
The success of the RB approach in parametrized boundary value problems
is only understood precisely from a mathematical perspective in a few very simple cases. This should motivate further theoretical investigations.

In this section, we discuss various tracks for the development of reduced basis techniques, both from a methodological viewpoint and in terms of possible applications, with a focus on the stochatic context presented above. 

% As a matter of fact, as we already said before,
% the RB method does not compete with the usual discretization techniques,
% like spectral or finite element methods,
% but rather tries to build upon optimal and robust approximations
% that are {\it a posteriori} adapted to
% complicated parametrized contexts already discretized.
% So there may just be little general theory to be expected.
% There is nevertheless a strong need for mathematical theory in the RB practice. This
%  in turn can influence RB methodological choices and then
% provide interesting specific computational reductions in a specific
% problem.
% We now discuss the development of such mathematical knowledge 
% with a focus on the stochastic applications.

\subsection{{\em A posteriori} estimation in the stochastic context}

We already emphasized that a crucial ingredient in the RB approach is
an accurate and fast {\it a posteriori estimator} for the approximation error
between two levels of discretization (the initially discretized, non-reduced one 
and the reduced one).
Therefore, before everything, the future developments of the RB method
should definitely concentrate on improving the {\em a posteriori} estimators.
In particular, for the new contexts of application that are {\it stochastic},
there seems to remain some room for a yet better understanding of the
{\em a posteriori} error estimation.
More precisely, the best way to evaluate the reduction error
when the Galerkin approximations (used by deterministic applications)
are replaced with Monte-Carlo approximations is still unclear.
For a first application of the RB ideas to stochastic applications,
we have used confidence intervals as 
a probabilistic measure of the Monte-Carlo approximation error. These confidence intervals are only reliable in the limit of infinitely many realizations of the random variables. But there are other possibilities, like using non-asymptotic upper-bounds for the error which hold whatever the number of realizations (using for example Chebyshev inequalities
or Berry-Esseen type bounds).
In addition, the numerical evaluation of the variance is not obvious either.
Until now, we have used Monte-Carlo estimators, but there exist other possibilities too
which could be faster or more accurate and should thus be tested. Finally, another idea related to the method presented in Section~\ref{sec:SDE}  would be to mimick the usual RB approach, by considering that the reference result is the one obtained with $M_{\rm large}$ realizations, and to develop {\em a posteriori} error bounds with respect to this reference solution (using for example conditional expectations with respect to the $M_{\rm large}$  realizations).

\subsection{Affine decompositions and the stochastic context}

As explained above, the RB approach is to date only efficient 
at yielding computational reductions in the context of {\em affine}
parametrization. However, as shown in the previous sections, a many-query parametrized framework is not necessarily
parametrized in an affine way.
So one may have to pretreat the problem in order to transform it
as the limit of a sequence of affinely-parametrized problems. It would thus be interesting to derive rapidly convergent affine approximations for non-affine problems.
For instance, the Karhunen-Loeve decomposition used to pretreat a random field
entering a partial differential equation as a coefficient may converge too slowly for an efficient use
of the RB method applied to truncated decompositions,
in the context of random fields with small correlation lengths.
One should then look for other possible affine representation of the random variations
in the input coefficient.
Now, there are many possible tracks to solve this problem,
like projecting the random field on a well-chosen basis 
for the oscillation modes of the coefficient for instance
(many possible bases may exist for the realizations of the
random field, depending on its regularity),
interpolating (recall the so-called {\it empirical interpolation} with magic points), etc.

\subsection{Application to Bayesian statistics}

A context where the RB ideas could be applied 
is {\it Bayesian statistics} where the many-query parametrized framework is naturally encountered. We now make this more precise by presenting a specific example that would be well-suited for the application of Algorithm 1 in Section~\ref{sec:SDE}.

Let us consider, for given values of a parameter $\mu_1$, an ensemble of observations $(x^{\mu_1}_i)_{1 \le i \le N^{\mu_1}_{\rm data}}$. Following the Bayesian framework, a stochastic model is proposed to model the observation: the quantities $(x^{\mu_1}_i)_{1 \le i \le N^{\mu_1}_{\rm data}}$ are supposed to form a set of independent and identically distributed samples following a given distribution parametrized by another set of parameters $\mu_2$ (think for example of a mixture of Gaussians, $\mu_2$ being then the triplets of weights, means and variances of each Gaussians). The Bayesian approach then consists in postulating a so-called {\em prior distribution} (with a probability density function denoted ${\rm Prior} (\mu_2)$ below) on the parameters $\mu_2$, and to compute the so-called {\em posterior distribution}, namely the distribution of $\mu_2$ given the observations (with a probability density function denoted $\Pi\big(\mu_2 | (x^{\mu_1}_i)_{1 \le i \le N^{\mu_1}_{\rm data}} \big)$ below). Of course, the posterior distribution for $\mu_2$ is expected to depend on $\mu_1$: for each $\mu_1$, the aim is thus to sample the probability measure~$\Pi\big(\mu_2 | (x^{\mu_1}_i)_{1 \le i \le N^{\mu_1}_{\rm data}} \big) \, d\mu_2$, with
$$\Pi\big(\mu_2 | (x^{\mu_1}_i)_{1 \le i \le N^{\mu_1}_{\rm data}} \big) = (Z^{\mu_1})^{-1} \Pi\big( (x^{\mu_1}_i)_{1 \le i \le N^{\mu_1}_{\rm data}} | \mu_2\big) \, {\rm Prior} (\mu_2) $$
where $\displaystyle Z^{\mu_1}$ is the normalization constant, and $\Pi\big( (x^{\mu_1}_i)_{1 \le i \le N^{\mu_1}_{\rm data}} | \mu_2\big)$ is the so-called likelihood function, namely the probability density function of the observations given the datas.
One possible technique to sample the posterior distribution consists in drawing samples according to the prior distribution, and to weight each of them using the likelihood function, which depends on $\mu_1$. With such a sampling technique, it is easy to draw correlated samples for various values of $\mu_1$. Following Algorithm 1 in Section~\ref{sec:SDE}, it would thus be possible to build a reduced basis based on the sampling of the posterior distribution for some selected values of $\mu_1$ (offline stage), in order to reduce the variance for the sampling of the posterior distribution for other generic values of $\mu_1$ (online stage).

\subsection{Relation to functional quantization}

One computationally demanding stochastic context that defines a many-query framework
is the approximation of the solution to a parametrized stochastic differential equation,
for many values of the parameter. We already mentioned applications in finance and rheology in Section 4, where a variance reduction technique based on RB was proposed. Another idea consists in first computing precise discretizations of a few processes
at some well-chosen parameter values, and then to use them for a faster
computation of an approximation of the processes for other values of the parameter.

{\em Functional quantization} is an approach that has independently been developed along this line, see for instance~\cite{luschy-pages-02,pages-printems-05}.
The idea of quantization is to approximate a square-integrable random variable with values in a Hilbert space by a random variable that takes a finite number of values, in an optimal way. In its simplest form, quantization deals with Gaussian random variables with values in $\R^d$, but it can also be applied to Gaussian processes.
The numerical approach developped in~\cite{pages-sellami-08} to solve stochastic differential equations is to first quantize the Brownian motion,
and then to solve a collection of ordinary differential equations in order to recover
approximations of the solutions to the stochastic differential equations
as linear combinations of the ordinary differential equations solutions.
Clearly, this approach for the discretization of stochastic differential equations
has intimate connection with a RB approach. In particular,
the computations are split into two parts:
an offline step, which is computationally expensive, to quantize the Brownian motion,
and then an online step to solve ordinary differential equations rather than stochastic differential equations.

In a setting where the stochastic differential equations are parametrized, a natural similar idea would be to quantize the solution to the stochastic differential equations for a few values of the parameter, and next to build the solution to the stochastic differential equation for another value of the parameter as a linear combination of these precomputed solutions.

% \bibliography{RB_book}

\begin{thebibliography}{OvdBH97}

\bibitem[AO00]{ainsworth00:_poster_error_estim_finit_elemen_analy}
M.~Ainsworth and J.~T. Oden.
\newblock {\em {\em A Posteriori\/} Error Estimation in Finite Element
  Analysis}.
\newblock Wiley-Interscience, 2000.


\bibitem[ASB78]{almroth78:_autom}
B.~O. Almroth, P.~Stern, and F.~A. Brogan.
\newblock Automatic choice of global shape functions in structural analysis.
\newblock {\em AIAA Journal}, 16:525--528, 1978.


\bibitem[Aro04]{arouna-03}
B.~Arouna.
\newblock Robbins-monroe algorithms and variance reduction in finance.
\newblock {\em The Journal of Computational Finance}, 7(2):35--62, 2004.


\bibitem[BNT07]{babuska-nobile-tempone-07}
I.~Babu\v{s}ka, F.~Nobile, and R.~Tempone.
\newblock A stochastic collocation method for elliptic partial differential
  equations with random input data.
\newblock {\em SIAM Journal on Numerical Analysis}, 45(3):1005--1034, 2007.


\bibitem[BTZ05]{babuska-tempone-zouraris-05}
I.~Babu\v{s}ka, R.~Tempone, and G.~Zouraris.
\newblock Solving elliptic boundary value problems with uncertain coefficients
  by the finite element method: the stochastic formulation.
\newblock {\em Comput. Meth. Appl. Mech. Eng.}, 194:1251--1294, 2005.


\bibitem[BNMP04]{barrault04:_empir_inter_method}
M.~Barrault, N.~C. Nguyen, Y.~Maday, and A.~T. Patera.
\newblock An ``empirical interpolation'' method: Application to efficient
  reduced-basis discretization of partial differential equations.
\newblock {\em C. R. Acad. Sci. Paris, S{\'e}rie I.}, 339:667--672, 2004.


\bibitem[BR95]{barrett95:_reduc_basis_method}
A.~Barrett and G.~Reddien.
\newblock On the reduced basis method.
\newblock {\em Z. Angew. Math. Mech.}, 75(7):543--549, 1995.

\bibitem[BR01]{becker-rannacher-01}
R.~Becker and R.~Rannacher.
\newblock An optimal control approach to a posteriori error estimation in
  finite element methods.
\newblock {\em Acta Numerica}, 37:1--225, 2001.


\bibitem[BP99]{bonvin-picasso-99}
J.~Bonvin and M.~Picasso.
\newblock Variance reduction methods for {CONNFFESSIT}-like simulations.
\newblock {\em J. Non-Newtonian Fluid Mech.}, 84:191--215, 1999.


\bibitem[Boy08]{boyaval-08}
S.~Boyaval.
\newblock Reduced-basis approach for homogenization beyond the periodic
  setting.
\newblock {\em SIAM Multiscale Modeling \& Simulation}, 7(1):466--494, 2008.

\bibitem[Boy09]{boyaval-these}
S.~Boyaval.
\newblock {\em Mathematical modeling and simulation for material science}.
\newblock PhD thesis, Universit{\'e} Paris Est, 2009.


\bibitem[BBM{\etalchar{+}}09]{boyaval-lebris-maday-nguyen-patera-09}
S.~Boyaval, C.~Le Bris, Y.~Maday, N.C. Nguyen, and A.T. Patera.
\newblock A reduced basis approach for variational problems with
stochastic
  parameters: Application to heat conduction with variable Robin
coefficient.
\newblock {\em Computer Methods in Applied Mechanics and Engineering},
  198(41--44):3187--3206, 2009.

\bibitem[BL09]{boyaval-lelievre-09}
{S}. {B}oyaval and {T}. {L}eli{\`e}vre.
\newblock A variance reduction method for parametrized stochastic differential
  equations using the reduced basis paradigm.
\newblock In Pingwen Zhang, editor, {\em Accepted for publication in
  Communication in Mathematical Sciences}, volume Special Issue ``Mathematical
  Issues on Complex Fluids'', 2009.

\bibitem[BMPPT10]{BMPPT10}
A. Buffa, Y. Maday, A.T. Patera, C. Prud'homme and G. Turinici.
\newblock A priori convergence of the greedy algorithm for the parametrized reduced basis.
\newblock Submitted to {\em M2AN (Math. Model. Numer. Anal.)}.

\bibitem[BGL06]{gunzburger06b}
J.~Burkardt, M.~D. Gunzburger, and H.~C. Lee.
\newblock {POD} and {CVT}-based reduced order modeling of {N}avier-{S}tokes
  flows.
\newblock {\em Comp. Meth. Applied Mech.}, 196:337--355, 2006.

\bibitem[Cia78]{Ciarlet}
Ph.\ Ciarlet.
\newblock {\em The Finite Element Method for Elliptic Problems}.
\newblock North-Holland, Amsterdam, 1978.

\bibitem[CDD]{cohen-perso}
A.~Cohen, R.~A. Devore, and W.~Dahmen.
\newblock Personal communication.

\bibitem[DBO01]{Deb01}
M.K. Deb, I.M. Babu{\v{s}}ka, and J.T. Oden.
\newblock Solution of stochastic partial differential equations using galerkin
  finite element techniques.
\newblock {\em Comp. Meth. Appl. Mech. Engrg.}, 190(48):6359--6372, 2001.

\bibitem[DNP{\etalchar{+}}04]{Debusschere04}
B.J. Debusschere, H.N. Najm, P.P. Pebay, O.M. Knio, R.G. Ghanem, and O.P.~Le
  Ma\^itre.
\newblock Numerical challenges in the use of polynomial chaos representations
  for stochastic processes.
\newblock {\em SIAM Journal on Scientific Computing}, 26(2):698--719, 2004.


\bibitem[Dep09]{deparis07}
S.~Deparis.
\newblock Reduced basis error bound computation of parameter-dependent
  {N}avier-{S}tokes equations by the natural norm approach.
\newblock {\em SIAM Journal of Numerical Analysis}, 46(4):2039--2067, 2009.

\bibitem[Dev93]{devore-93}
R.~A. Devore.
\newblock Constructive approximation.
\newblock {\em Acta Numerica}, 7:51--150, 1993.


\bibitem[FR83]{fink83}
J.~P. Fink and W.~C. Rheinboldt.
\newblock On the error behavior of the reduced basis technique for nonlinear
  finite element approximations.
\newblock {\em Z. Angew. Math. Mech.}, 63(1):21--28, 1983.

\bibitem[FM71]{fox71:_approx_analy_techn_desig_calcu}
R.~L. Fox and H.~Miura.
\newblock An approximate analysis technique for design calculations.
\newblock {\em AIAA Journal}, 9(1):177--179, 1971.


\bibitem[GS91]{ghanem-spanos-03}
R.~G. Ghanem and P.~D. Spanos.
\newblock {\em Stochastic Finite Elements: A Spectral Approach}.
\newblock Springer {V}erlag, {N}ew {Y}ork edition, 1991.
\newblock Revised {F}irst {D}over {E}dition, 2003.


\bibitem[GMNP07]{m2an_magic}
M.~A. Grepl, Y.~Maday, N.~C. Nguyen, and A.~T. Patera.
\newblock Efficient reduced-basis treatment of nonaffine and nonlinear partial
  differential equations.
\newblock {\em M2AN (Math. Model. Numer. Anal.)}, 41(2):575--605, 2007.
\newblock (doi: 10.1051/m2an:2007031).

\bibitem[GP05]{grepl04:_reduc_basis_approx_time_depen}
M.~A. Grepl and A.~T. Patera.
\newblock {\em A Posteriori} error bounds for reduced-basis approximations of
  parametrized parabolic partial differential equations.
\newblock {\em M2AN (Math. Model. Numer. Anal.)}, 39(1):157--181, 2005.

\bibitem[GvL96]{golub-vanloan-96}
G.~Golub and C.~van Loan.
\newblock {\em Matrix computations, third edition}.
\newblock The Johns Hopkins University Press, London, 1996.

\bibitem[Gun89]{gunzburger89:_finit}
M.~D. Gunzburger.
\newblock {\em Finite Element Methods for Viscous Incompressible Flows}.
\newblock Academic Press, 1989.

\bibitem[HO08a]{Haasdonk-08}
B.~Haasdonk and M.~Ohlberger.
\newblock Reduced basis method for finite volume approximations of parametrized
  linear evolution equations.
\newblock {\em M2AN (Math. Model. Numer. Anal.)},
  42(3):277--302, 2008.

\bibitem[HO08b]{Haasdonk-Ohlberger-08-2}
B.~Haasdonk and M.~Ohlberger.
\newblock Adaptive Basis Enrichment for the Reduced Basis Method
Applied to Finite Volume Schemes.
\newblock  In Proc. 5th International Symposium on Finite Volumes
for Complex Applications, June 08-13, 2008, Aussois, France.
\newblock 471--478, 2008.

\bibitem[HDH64]{hammersley-handscomb-64}
J.~Hammersley and eds. D.~Handscomb.
\newblock {\em Monte Carlo Methods}.
\newblock Chapman and Hall Ltd, London, 1964.

\bibitem[HKY{\etalchar{+}}09]{huyhn-knezevic-chen-hesthaven-patera-09}
D.B.P. Huynh, D.J. Knezevic, Y.Chen, J.S. Hesthaven, and A.~T. Patera.
\newblock A natural-norm successive constraint method for inf-sup lower bounds.
\newblock {\em Computer Methods in Applied Mechanics and Engineering}, 2009.
\newblock submitted, preprint version available on URL
  {\tt http://augustine.mit.edu/methodology/papers/}.

\bibitem[HRSP07]{huynh-rozza-sen-patera-07}
D.B.P. Huynh, G.~Rozza, S.~Sen, and A.T. Patera.
\newblock A successive constraint linear optimization method for lower bounds
  of parametric coercivity and inf-sup stability constants.
\newblock {\em C. R. Math. Acad. Sci. Paris}, 345:473--478, 2007.


\bibitem[HP07]{Patera_Huynh06}
D.~B.~P. Huynh and A.~T. Patera.
\newblock Reduced-basis approximation and {\em a posteriori\/} error estimation
  for stress intensity factors.
\newblock {\em Int. J. Num. Meth. Eng.}, 72(10):1219--1259, 2007.
\newblock (doi: 10.1002/nme.2090).

\bibitem[Kar46]{k.46:_zur_spekt_prozes}
K.~Karhunen.
\newblock Zur spektraltheorie stochastischer prozesse.
\newblock {\em Annales Academiae Scientiarum Fennicae}, 37, 1946.

\bibitem[KP00]{kloeden-platen-00}
P.~Kloeden and E.~Platen.
\newblock {\em Numerical Solution of Stochastic Differential Equations}.
\newblock Springer, 2000.

\bibitem[KP09]{knezevic-patera-09}
D.J. Knezevic and A.T. Patera.
\newblock A certified reduced basis method for the Fokker-Planck equation of
  dilute polymeric fluids: FENE dumbbells in extensional flow.
\newblock submitted to SIAM Journal of Scientific Computing, 2009.

\bibitem[KV02]{kunish02}
K.~Kunisch and S.~Volkwein.
\newblock Galerkin proper orthogonal decomposition methods for a general
  equation in fluid dynamics.
\newblock {\em SIAM J. Num. Analysis}, 40(2):492--515, 2002.

\bibitem[LB05]{le_bris-05}
C. Le~Bris.
\newblock {\em Syst\`emes multi-\'echelles : Mod\'elisation \& simulation},
  volume~47 of {\em Math\'ematiques \& applications}.
\newblock Springer, 2005.

\bibitem[LL07]{LL07}
C. Le Bris and T. Leli\`evre.
\newblock Multiscale modelling of complex fluids: A mathematical initiation.
\newblock In Multiscale Modeling and Simulation in Science Series, B. Engquist, P. Lötstedt, O. Runborg, eds.,
\newblock LNCSE 66, Springer, p. 49-138, 2009.

% \bibitem[LBLM08]{BRIS:2008:INRIA-00336911:1}
% {C}. {L}e {B}ris, {T}. {L}eli{\`e}vre, and {Y}. {M}aday.
% \newblock {R}esults and questions on a nonlinear approximation approach for
%   solving high-dimensional partial differential equations.
% \newblock {\em (preprint submitted for publication
%   http://hal.inria.fr/inria-00336911/en/)}, 2008.

\bibitem[LL02]{lienhard}
J.~H. {Lienhard IV} and J.~H. {Lienhard V}.
\newblock {\em A {H}eat {T}ransfer {T}extbook}.
\newblock {P}hlogiston {P}ress, {C}ambridge, {M}ass., 2002.


\bibitem[Lin91]{lee91:_estim}
M.~Y.~Lin Lee.
\newblock Estimation of the error in the reduced-basis method solution of
  differential algebraic equations.
\newblock {\em SIAM Journal of Numerical Analysis}, 28:512--528, 1991.

\bibitem[Lo{\`e}78]{loeve-78}
M.~Lo{\`e}ve.
\newblock {\em Probability Theory}, volume I-II.
\newblock Springer, {N}ew {Y}ork, 1978.

\bibitem[LuPa02]{luschy-pages-02}
H.~Luschy and G.~Pag{\`e}s
\newblock Functional quantization of Gaussian processes
\newblock {\em J. Funct. Anal.}, 196:486--531, 2002.

\bibitem[Mad06]{maday06}
Y.~Maday.
\newblock Reduced--basis method for the rapid and reliable solution of partial
  differential equations.
\newblock In {\em Proceedings of International Conference of Mathematicians,
  Madrid}. European Mathematical Society Eds., 2006.

\bibitem[MNPP09]{maday-nguyen-patera-pau-07}
Y.~Maday, N.~C. Nguyen, A.~T. Patera, and G.~Pau.
\newblock A general, multipurpose interpolation procedure: the magic points.
\newblock {\em Communications on Pure and Applied Analysis}, 8(1):383--404,
  2009.

\bibitem[MPT02a]{maday02:_prior_conver_theor_reduc}
Y.~Maday, A.~T. Patera, and G.~Turinici.
\newblock {\em A Priori\/} convergence theory for reduced-basis approximations
  of single-parameter elliptic partial differential equations.
\newblock {\em Journal of Scientific Computing}, 17(1-4):437--446, 2002.

\bibitem[MPT02b]{maday-02}
Y.~Maday, A.T. Patera, and G.~Turinici.
\newblock Global a priori convergence theory for reduced-basis approximations
  of single-parameter symmetric coercive elliptic partial differential
  equations.
\newblock {\em C. R. Acad. Sci. Paris, Ser. I}, 335(3):289--294, 2002.

\bibitem[MK05]{matthies-keese-04}
H.G. Matthies and A.~Keese.
\newblock Galerkin methods for linear and nonlinear elliptic stochastic pdes.
\newblock {\em Comput. Meth. Appl. Mech. Eng.}, 194:1295--1331, 2005.


\bibitem[MO95]{melchior-ottinger-95}
M.~Melchior and H.C. \"Ottinger.
\newblock Variance reduced simulations of stochastic differential equations.
\newblock {\em J. Chem. Phys.}, 103:9506--9509, 1995.


\bibitem[MT06]{milstein-tretyakov-06}
G.N. Milstein and M.V. Tretyakov.
\newblock Practical variance reduction via regression for simulating
  diffusions.
\newblock Technical Report MA-06-19, School of Mathematics and Computer
  Science, University of Leicester, 2006.


\bibitem[Nag79]{nagy79:_modal}
D.~A. Nagy.
\newblock Modal representation of geometrically nonlinear behaviour by the
  finite element method.
\newblock {\em Computers and Structures}, 10:683--688, 1979.


\bibitem[Nguyen07]{nguyen-07}
N.~C. Nguyen.
\newblock A multiscale reduced-basis method for parametrized elliptic
partial differential equations with multiple scales.
\newblock {\em Journal of Computational Physics}, 227:9807--9822, 2007.

\bibitem[Noo81]{noor81:_recen}
A.~K. Noor.
\newblock Recent advances in reduction methods for nonlinear problems.
\newblock {\em Comput. Struct.}, 13:31--44, 1981.

\bibitem[Noo82]{noor82}
A.~K. Noor.
\newblock On making large nonlinear problems small.
\newblock {\em Comp. Meth. Appl. Mech. Engrg.}, 34:955--985, 1982.


\bibitem[Noor84]{noor84:_reduc}
A.~K. Noor, C.~D. Balch, and M.~A. Shibut.
\newblock Reduction methods for non-linear steady-state thermal analysis.
\newblock {\em Int. J. Num. Meth. Engrg.}, 20:1323--1348, 1984.


\bibitem[NP80]{noor80:_reduc}
A.~K. Noor and J.~M. Peters.
\newblock Reduced basis technique for nonlinear analysis of structures.
\newblock {\em AIAA Journal}, 18(4):455--462, 1980.


\bibitem[NP83a]{noor83:_multip}
A.~K. Noor and J.~M. Peters.
\newblock Multiple-parameter reduced basis technique for bifurcation and
  post-buckling analysis of composite plates.
\newblock {\em Int. J. Num. Meth. Engrg.}, 19:1783--1803, 1983.

\bibitem[NP83b]{noor83:_recen}
A.~K. Noor and J.~M. Peters.
\newblock Recent advances in reduction methods for instability analysis of
  structures.
\newblock {\em Comput. Struct.}, 16:67--80, 1983.

\bibitem[NPA84]{noor84:_mixed}
A.~K. Noor, J.~M. Peters, and C.~M. Andersen.
\newblock Mixed models and reduction techniques for large-rotation nonlinear
  problems.
\newblock {\em Comp. Meth. Appl. Mech. Engrg.}, 44:67--89, 1984.


\bibitem[New96]{newman96}
A.~J. Newman.
\newblock Model reduction via the {K}arhunen-{L}oeve expansion part i: an
  exposition.
\newblock {\em Technical Report Institute for System Research University of
  Maryland}, (96-322), 1996.

\bibitem[New94]{newton-94}
N.~J. Newton.
\newblock Variance reduction for simulated diffusions.
\newblock {\em SIAM J. Appl. Math.}, 54(6):1780--1805, 1994.

\bibitem[NRHP09]{nguyen-09}
N.~C. Nguyen, G.~Rozza, D.~B.~P. Huynh, and A.~T Patera.
\newblock Reduced basis approximation and a posteriori error estimation for
  parametrized parabolic pdes; application to real-time bayesian parameter
  estimation.
\newblock In Biegler, Biros, Ghattas, Heinkenschloss, Keyes, Mallick, Tenorio,
  van Bloemen~Waanders, and Willcox, editors, {\em Computational Methods for
  Large Scale Inverse Problems and Uncertainty Quantification}, John Wiley \&
  Sons, UK, 2009.

\bibitem[NRP08]{Calcolo}
N.~C. Nguyen, G.~Rozza, and A.~T. Patera.
\newblock Reduced basis approximation and a posteriori error estimation for the
  time-dependent viscous burgers equation.
\newblock {\em Calcolo}, 46(3):157--185, 2009.

\bibitem[NVP05]{nguyen-veroy-patera-05}
N.~C. Nguyen, K.~Veroy, and A.T. Patera.
\newblock {\em Certified real-time solution of parametrized partial
  differential equations}, pages 1523--1558.
\newblock {S}pringer, 2005.
\newblock (S. Yip, editor).

\bibitem[Nou07]{nouy-07-1}
A.~Nouy.
\newblock A generalized spectral decomposition technique to solve a class of
  linear stochastic partial differential equations.
\newblock {\em Computer Methods in Applied Mechanics and Engineering},
  196(45-48):4521--4537, 2007.

\bibitem[Nou08]{nouy-07-2}
A.~Nouy.
\newblock Generalized spectral decomposition method for solving stochastic
  finite element equations: Invariant subspace problem and dedicated
  algorithms.
\newblock {\em Computer Methods in Applied Mechanics and Engineering},
  197(51-52):4718--4736, 2008.


\bibitem[OvdBH97]{ottinger-vandenbrule-hulsen-97}
H.~C. \"{O}ttinger, B.~H. A.~A. van~den Brule, and M.~A. Hulsen.
\newblock Brownian configuration fields and variance reduced {CONNFFESSIT}.
\newblock {\em J. Non-Newtonian Fluid Mech.}, 70(30):255 -- 261, 1997.

\bibitem[PP05]{pages-printems-05}
G.~Pag{\`e}s and J.~Printems.
\newblock Functional quantization for numerics with an application to option
  pricing.
\newblock {\em Monte Carlo Methods and Appl.}, 11(4):407--446, 2005.

\bibitem[PS08]{pages-sellami-08}
G.~Pag{\`e}s and A.~Sellami.
\newblock Convergence of multi-dimensional quantized SDE's.
\newblock {\tt http://arxiv.org/abs/0801.0726}, 2008.

\bibitem[PR07a]{Patera_Ronquist07}
A.~T. Patera and E.~M. R{\o}nquist.
\newblock Reduced basis approximations and {\em a posteriori\/} error
  estimation for a {B}oltzmann model.
\newblock {\em Computer Methods in Applied Mechanics and Engineering},
  196:2925--2942, 2007.

\bibitem[PR07b]{patera_rozza}
A.~T. Patera and G.~Rozza.
\newblock {\em Reduced Basis Approximation and A Posteriori Error Estimation
  for Parametrized Partial Differential Equations}.
\newblock Copyright MIT, 2006--2007.
\newblock To appear in MIT Pappalardo Monographs in Mechanical Engineering.


\bibitem[Pet89]{peterson89}
J.~S. Peterson.
\newblock The reduced basis method for incompressible viscous flow
  calculations.
\newblock {\em SIAM J. Sci. Stat. Comput.}, 10(4):777--786, 1989.

\bibitem[Pin85]{pinkus}
A.~Pinkus.
\newblock {\em n-{W}idths in Approximation Theory}.
\newblock Springer, 1985.


\bibitem[Por85]{porsching85:_estim}
T.~A. Porsching.
\newblock Estimation of the error in the reduced basis method solution of
  nonlinear equations.
\newblock {\em Mathematics of Computation}, 45(172):487--496, 1985.

\bibitem[Por87]{porsching87}
T.~A. Porsching and M.~Y.~Lin Lee.
\newblock The reduced-basis method for initial value problems.
\newblock {\em SIAM Journal of Numerical Analysis}, 24:1277--1287, 1987.


\bibitem[PRV{\etalchar{+}}02]{prud'homme02:_reliab_real_time_solut_param}
C.~Prud'homme, D.~Rovas, K.~Veroy, Y.~Maday, A.~T. Patera, and G.~Turinici.
\newblock Reliable real-time solution of parametrized partial differential
  equations: Reduced-basis output bounds methods.
\newblock {\em Journal of Fluids Engineering}, 124(1):70--80, 2002.

\bibitem[Qua09]{quarteroni-09}
A.~Quarteroni.
\newblock {\em Numerical Models for Differential Problems}, volume~2.
\newblock Springer, 2009.


\bibitem[Rhe81]{rheinboldt81:_numer_analy_contin_method_nonlin_struc_probl}
W.~C. Rheinboldt.
\newblock Numerical analysis of continuation methods for nonlinear structural
  problems.
\newblock {\em Computers and Structures}, 13(1-3):103--113, 1981.

\bibitem[Rhe93]{rheinboldt93:_theor_error_estim_reduc_basis}
W.~C. Rheinboldt.
\newblock On the theory and error estimation of the reduced basis method for
  multi-parameter problems.
\newblock {\em Nonlinear Analysis, Theory, Methods and Applications},
  21(11):849--858, 1993.


\bibitem[RHP]{Rozza08:arcme}
G.~Rozza, D.B.P. Huynh, and A.~T. Patera.
\newblock Reduced basis approximation and a posteriori error estimation for
  affinely parametrized elliptic coercive partial differential equations ---
  application to transport and continuum mechanics.
\newblock {\em Arch. Comput. Methods Eng.}
\newblock to appear September 2008.


\bibitem[ST06]{schwab-todor-06}
C.~Schwab and R.A. Todor.
\newblock Karhunen-lo\`{e}ve approximation of random fields by generalized fast
  multipole methods.
\newblock {\em Journal of Computational Physics}, 217(1):100--122, 2006.

\bibitem[Sen08]{sen-08}
S.~Sen.
\newblock Reduced-basis approximation and a posteriori error estimation for
  many-parameter heat conduction problems.
\newblock {\em Numerical Heat Transfer, Part B: Fundamentals}, 54(5), 2008.


\bibitem[SVH{\etalchar{+}}06]{Pat06:naturalNorms}
S.~Sen, K.~Veroy, D.~B.~P. Huynh, S.~Deparis, N.~C. Nguyen, and A.~T. Patera.
\newblock ``{N}atural norm'' {\em a posteriori\/} error estimators for reduced
  basis approximations.
\newblock {\em Journal of Computational Physics}, 217:37--62, 2006.

\bibitem[SF73]{strang_fix}
G.~Strang and G.~J. Fix.
\newblock {\em An Analysis of the Finite Element Method}.
\newblock Prentice-Hall, 1973.

\bibitem[V.N08]{temlyakov-08}
V.N.Temlyakov.
\newblock Nonlinear methods of approximation.
\newblock {\em Foundations of Computational Mathematics}, 3:33--107, 2008.

\bibitem[VPRP03]{veroy03:_poster_error_bound_reduc_basis}
K.~Veroy, C.~Prud'homme, D.~V. Rovas, and A.~T. Patera.
\newblock {\em A Posteriori\/} error bounds for reduced-basis approximation of
  parametrized noncoercive and nonlinear elliptic partial differential
  equations.
\newblock In {\em Proceedings of the 16th AIAA Computational Fluid Dynamics
  Conference}, 2003.
\newblock Paper 2003-3847.

\bibitem[VRP02]{veroy_lions}
K.~Veroy, D.~Rovas, and A.~T. Patera.
\newblock {\em A Posteriori\/} error estimation for reduced-basis approximation
  of parametrized elliptic coercive partial differential equations: ``convex
  inverse'' bound conditioners.
\newblock {\em ESAIM: Control, Optimization and Calculus of Variations},
  8:1007--1028, 2002.

\end{thebibliography}
% \bibliographystyle{plain}

\newcommand{\etalchar}[1]{$^{#1}$}

\end{document}